\documentclass{article}
\usepackage{amssymb}

\newcommand{\PREP}[1]{#1}
\newcommand{\WSCI}[1]{}

\newtheorem{thm}{Theorem}[section]
\newtheorem{prop}[thm]{Proposition}
\newtheorem{lem}[thm]{Lemma}
\newtheorem{cor}[thm]{Corollary}

\newcommand{\rem}[1]{
\par\medskip\noindent{\bf Remark: }\ #1
}

\PREP{
\newcommand{\eqref}[1]{({#1})}
}
\WSCI{
\newcommand{\eqref}[1]{Eq.\,({#1})}
}

\newcommand{\BZ}{{\mathbb Z}\,}
\newcommand{\BQ}{{\mathbb Q}}
\newcommand{\BC}{{\mathbb C}}
\newcommand{\CA}{{\mathcal A}}
\newcommand{\CB}{{\mathcal B}}

\newcommand{\CK}{{\mathcal K}}

\newcommand{\CO}{{\mathcal O}}
\newcommand{\CR}{{\mathcal R}}
\newcommand{\CS}{{\mathcal S}}

\newcommand{\frg}{{\mathfrak g}}
\newcommand{\frh}{{\mathfrak h}}
\newcommand{\frn}{{\mathfrak n}}
\newcommand{\frb}{{\mathfrak b}}

\newcommand{\frS}{{\mathfrak S}}
\newcommand{\ds}{\dot{s}}

\newcommand{\br}[1]{{\langle{#1}\rangle}}
\newcommand{\ad}{\mbox{\rm ad}}
\newcommand{\Ad}{\mbox{\rm Ad}}
\newcommand{\pbr}[1]{{\{{#1}\}}}
\newcommand{\pad}{{\mbox{\rm ad}}_{\{\}}}
\newcommand{\iso}{\stackrel{\sim}{\to}}
\newcommand{\comment}[1]{}

\newcommand{\hW}{{\mathcal W}}
\newcommand{\Aut}{\mbox{\rm Aut}}
\newcommand{\Hom}{\mbox{\rm Hom}}
\newcommand{\End}{\mbox{\rm End}}
\newcommand{\Spec}{\mbox{\rm Spec}}
\newcommand{\Sym}{\mbox{\rm Sym}}
\newcommand{\projlim}[1]{
\begin{picture}(18,8)(0,0)
\put(0,0){$\lim$}
\put(0,-5){$\longleftarrow$}
\put(15,-5){\scriptsize $#1$}
\end{picture}
\,}
\newcommand{\dfrac}[2]{{\displaystyle\frac{#1}{#2}}}
\newcommand{\QED}{{\hfill Q.E.D.}}

\begin{document}
\PREP{
\title{\bf
Birational Weyl group action arising from 
a nilpotent Poisson algebra
}
\author{Masatoshi NOUMI and Yasuhiko YAMADA \\
{\normalsize Department of Mathematics, Kobe University}\\
{\normalsize Rokko, Kobe 657-8501, Japan}
} 
\date{}
\maketitle
\begin{abstract}
We propose a general method to realize an arbitrary 
Weyl group of Kac-Moody type as a group of birational canonical 
transformations, by means of a nilpotent Poisson algebra. 
\end{abstract}
}

\WSCI{
\title{
Birational Weyl group action arising from 
a nilpotent Poisson algebra
}
\author{Masatoshi NOUMI} 
\address{
Department of Mathematics, Kobe University, 
Rokko, Kobe 657-8501, Japan\\
E-mail: noumi@math.sci.kobe-u.ac.jp
}
\author{Yasuhiko YAMADA}
\address{
Department of Mathematics, Kobe University, 
Rokko, Kobe 657-8501, Japan\\
E-mail: yamaday@math.sci.kobe-u.ac.jp
}

\maketitle

\abstracts{
We propose a general method to realize an arbitrary 
Weyl group of Kac-Moody type as a group of birational canonical 
transformation, by means of a nilpotent Poisson algebra. 
}
}

\section*{Introduction}

In this paper we propose a general method to realize an arbitrary 
Weyl group of Kac-Moody type as a group of birational canonical 
transformations. 
Our construction is formulated by means of a nilpotent Poisson 
algebra. 
It can be regarded as a conceptual generalization of the 
birational Weyl group actions
proposed in our previous paper\WSCI{.\cite{NY} }\PREP{ \cite{NY}.} 
We also discuss a certain cocycle related to this realization
and its regularity, and give a proof to a generalization of 
the regularity conjecture\WSCI{.\cite{NY} }\PREP{ \cite{NY}.}
\par\medskip
The plan of this paper is as follows. 
We give a summary of our main results in Section 1. 
Fixing a generalized Cartan matrix $A$, 
we take as a datum a nilpotent Poisson algebra 
$\CA_0$ generated by a set of elements $\varphi_i$ ($i\in I$) 
satisfying the Serre relations (specified by the GCM $A$) with respect 
to the adjoint action by the Poisson bracket. 
Starting from such an $\CA_0$, 
we formulate a method 
to realize the Weyl group $W=W(A)$ associated with $A$,
as a group of birational canonical transformations 
of a field of rational functions defined by $\CA_0$ 
(Theorem \ref{thm-A}). 
We also introduce $\tau$-functions in Theorem \ref{thm-B} for 
our realization, and formulate in Theorem \ref{thm-C} a certain 
regularity property of the {\em $\tau$-cocycle} arising from 
the transformations of $\tau$-functions.
After describing explicit examples 
in the cases of rank 2, 
we give in Section 3 
a proof of Theorems \ref{thm-A} and \ref{thm-B}. 

In Section 4, we explain a Lie theoretic background 
of our birational realization of the Weyl group, 
in terms of Kac-Moody 
Lie algebras and Kac-Moody groups. 
In fact, we consider 
the birational {\em dressing} action of a lift $\dot{W}$ 
of the Weyl group 
on the Borel subgroup, induced through 
the Gauss decomposition in the Kac-Moody group. 
Our realization is then obtained by transferring this dressing  
action to the Borel subalgebra through the adjoint action. 
We give a proof of regularity of the $\tau$-cocycle 
(Theorem \ref{thm-C}) in Section 5 by using the 
geometric interpretation of Section 4. 
Finally in Section 6, we give some remarks related to 
our birational realization of Weyl groups. 

\section{Summary of results} 
\subsection{Birational realization of the Weyl group}
Let $A=(a_{ij})_{i,j\in I}$ be a 
generalized Cartan matrix (GCM for short). 
By definition, $A$ is an integer matrix satisfying the conditions
\begin{equation}
a_{jj}=2;\quad a_{ij}\le 0 \quad(i\ne j);\quad
a_{ij}=0\ \ \Longleftrightarrow\ \ a_{ji}=0,
\end{equation}
for any $i,j\in I$. 
When the indexing set is infinite, we always assume that 
$A$ is {\em locally finite\/}; i.e., for any $i\in I$, 
$a_{ij}=0$ except for a finite number of $j$'s. 
We denote 
the {\em root lattice} for $A$ 
by $Q=\bigoplus_{i\in I}\,\BZ\alpha_i$
and the {\em coroot lattice}
by $Q^\vee=\bigoplus_{i\in I}\BZ h_i$, 
where $\alpha_i$ and $h_i$ are the simple roots and the 
simple coroots, respectively. 
The canonical pairing $\br{\,,\,} :Q^\vee\times Q\to \BZ$ 
between the two lattices is defined 
by $\br{h_i,\alpha_j}=a_{ij}$  ($i,j\in I$).  
The Weyl group $W(A)$ for $A$ is defined by the generators 
$r_i$ ($i\in I$) with fundamental relations 
\begin{equation}
r_i^2=1,\qquad (r_i\,r_j)^{m_{ij}}=1\quad(i\ne j),
\end{equation}
where $m_{ij}=2,3,4,6$ or $\infty$ according as 
$a_{ij}a_{ji}=0,1,2,3,$ or $\ge 4$. 
This group acts naturally on 
$Q$ and $Q^\vee$ by 
\begin{equation}
r_i.\alpha=\alpha-\alpha_i\br{h_i,\alpha}\quad(\alpha\in Q),
\quad
r_i.h=h-\br{h,\alpha_i}h_i\quad(h\in Q^\vee),
\end{equation}
respectively, 
so that $\br{w.h,w.\alpha}=\br{h,\alpha}$ for any 
$h\in Q^\vee$ and $\alpha\in Q$. 
Let $\BC[\lambda]$ be the polynomial ring in the indeterminates 
$\lambda=(\lambda_i)_{i\in I}$.  
In the following context, each $\lambda_i$ will be regarded as a 
variable corresponding to the {simple coroot} $h_i$.
The Weyl group $W(A)=\br{\,r_i \,(i\in I)}$ acts on $\BC[\lambda]$ 
as a group of automorphisms such that 
$r_j(\lambda_i)=\lambda_i-a_{ij}\lambda_j$ for $i,j\in I$.

By a {\em Poisson algebra}, we mean a commutative 
$\BC$-algebra $\CA$ endowed with 
a skew-symmetric bilinear form
$\pbr{\,,\,}: \CA\times\CA\to\CA$,
called the {\em Poisson bracket},
such that
\begin{equation}\begin{array}{cl}
\smallskip
\mbox{(i)} &
\pbr{fg,h}=\pbr{f,h}g+f\pbr{g,h},\quad
\pbr{f,gh}=\pbr{f,g}h+g\pbr{f,h},\cr
\mbox{(ii)} &
\pbr{f,\pbr{g,h}}+\pbr{g,\pbr{h,f}}+\pbr{h,\pbr{f,g}}=0,
\end{array}\end{equation}
for any $f,g,h\in \CA$.
A homomorphism $T: \CA\to\CB$ between two Poisson algebras will 
be called a {\em canonical transformation\/}:
\begin{equation}
T(fg)=T(f)T(g),\quad T(\pbr{f,g})=\pbr{T(f),T(g)}\quad(f,g\in \CA).
\end{equation}
\par\medskip
We now fix a generalized Cartan matrix $A=(a_{ij})_{i,j\in I}$, and 
a Poisson algebra $\CA_0$. 
We assume that the algebra $\CA_0$ has no zerodivisors, and that
as a Poisson algebra, $\CA_0$ is generated by 
a set of nonzero elements $(\varphi_i)_{i\in I}$ such that 
\begin{equation}
\pad(\varphi_i)^{-a_{ij}+1}(\varphi_j)=0\quad(i\ne j),
\end{equation}
where $\pad(f)=\pbr{f,\cdot}$ stands for the adjoint action by 
the Poisson bracket. 
We denote by
$\CA=\CA_0[\lambda]$  the ring of polynomials
in the $\lambda$-variables with coefficients in $\CA_0$,  
and by $\CK=Q(\CA)$ the field of fractions of $\CA$.   
Roughly speaking, $\CK$ is 
the field of rational functions 
in the variables $\lambda_j$, and $\varphi_j$ together with 
the Poisson brackets among $\varphi_j$'s: 
\begin{equation}
\CK=\BC(
\lambda_j; 
\varphi_j, \pbr{\varphi_i,\varphi_j},
\pbr{\varphi_i,\pbr{\varphi_j,\varphi_k}}, \ldots).
\end{equation}
Note that the Poisson bracket $\pbr{\,,\,}$ 
of $\CA_0$ extends uniquely  to $\CK$  so that $\pbr{\lambda_i,\varphi}=0$ 
($\varphi\in\CK$)  for any $i\in I$. 
As to the Weyl group action on the $\lambda$-variables, 
we use the same notation $r_i$ for the $\CA_0$-linear automorphism of 
$\CK$ defined by $r_i(\lambda_j)=\lambda_j-a_{ji}\lambda_i$ ($j\in I$).
For each $i\in I$, we define a homomorphism 
$s_i :\CA\to\CA[\,\varphi_j^{-1} (j\in I)]$ as the composition
\begin{equation}
s_i=\exp\left(\frac{\lambda_i}{\varphi_i}\pad(\varphi_i)\right)\circ r_i
\quad(i\in I). 
\end{equation}
Note that, for any $\psi\in\CA_0$, the action of $s_i$ on $\psi$ is 
determined as the {\em finite} sum
\begin{equation}
s_i(\psi)=\psi+\frac{\lambda_i}{\varphi_i}\pbr{\varphi_i,\psi}+
\frac{1}{2!}\left(\frac{\lambda_i}{\varphi_i}\right)^2
\pbr{\varphi_i,\pbr{\varphi_i,\psi}}+\cdots,
\end{equation}
since the action of $\pad(\varphi_i)$ is locally nilpotent on $\CA_0$. 
These homomorphisms $s_i :\CA\to\CA[\,\varphi_j^{-1} (j\in I)]$ ($i\in I$) 
extend to automorphisms of $\CK=Q(\CA)$, for which 
we use the same notation $s_i$.

\WSCI{\par\medskip}
\begin{thm}\label{thm-A}
The automorphism $s_i$ $(i\in I)$ of the field of fractions 
$\CK=Q(\CA)$, defined as above, 
give a realization of the Weyl group $W(A)$ for the GCM $A$, 
as a group of canonical transformations of $\CK$. 
Namely, 
\begin{enumerate}
\item[$(1)$]
These automorphisms preserve the Poisson 
bracket of $\CK$ $:$  For each $i\in I$, 
\begin{equation}
s_i(\pbr{\varphi,\psi})=\pbr{s_i(\varphi), s_i(\psi)}
\qquad(\varphi,\psi\in \CK). 
\end{equation}
\item[$(2)$] 
They satisfy the fundamental relations for the generators of $W(A)$:
\begin{equation}
s_i^2=1,\qquad (s_i\,s_j)^{m_{ij}}=1\quad(i\ne j),
\end{equation}
where $m_{ij}=2,3,4,6$ or $\infty$ according as 
$a_{ij}a_{ji}=0,1,2,3,$ or $\ge 4$. 
\end{enumerate}
\end{thm}
\WSCI{\par\medskip}
\noindent
Theorem \ref{thm-A} is 
a systematic generalization of the realization of the Weyl group 
we discussed previously\WSCI{,\cite{NY} }\PREP{ \cite{NY},}
in terms of nilpotent 
Poisson algebras. 
(See Remark at the end of this section.)

\par\medskip
Our realization of the Weyl group 
is closely related to the universal exponential solution to 
the Yang-Baxter equation, due  to S.~Fomin and A.N.~Kirillov\WSCI{.\cite{FK} }
\PREP{ \cite{FK}. }
In order to clarify the point, let us define the {\em $R$-operator} 
$R_i(t)$ with a formal parameter $t$ by
\begin{equation}
R_i(t)=\exp\left(\frac{t}{\varphi_i}\pad(\varphi_i)\right)
\end{equation}
for each $i\in I$, so that $s_i=R_i(\lambda_i)\circ r_i$. 
The Coxeter relation $(s_i\,s_j)^{m_{ij}}=1$ ($i\ne j$) 
is then equivalent to the following {\em Yang-Baxter equation}, 
according to the type of the root system of rank two 
defined by $\alpha_i$ and $\alpha_j$:
\begin{equation}\label{YBE}
\begin{array}{cc}
\smallskip
(2A_1)&  R_i(a) R_j(b)= R_j(b) R_i(a),\cr
\smallskip
(A_2)&R_i(a) R_j(a+b) R_i(b)=R_j(b) R_i(a+b) R_j(a),\cr
(B_2)&R_i(a) R_j(a+b) R_i(a+2b) R_j(b)\cr
\smallskip
&=R_j(b) R_i(a+2b) R_j(a+b) R_i(a),\cr
(G_2)&R_i(a) R_j(a+b) R_i(2a+3b) R_j(a+2b) R_i(a+3b) R_j(b)\cr
&=R_j(b) R_i(a+3b) R_j(a+2b) R_i(2a+3b) R_j(a+b) R_i(a),
\end{array}
\end{equation}
where $a=\lambda_i$ and $b=\lambda_j$. 
These four cases corresponds to the values 
$(a_{ij},a_{ji})=(0,0), (-1,-1), (-2,-1), (-3,-1)$, respectively.

\subsection{$\tau$-Functions and the $\tau$-cocycle}
We now introduce a new set of variables $\tau_i$ ($i\in I$), 
called {\em $\tau$-functions}, 
and extend the action of the Weyl group to these variables. 
We denote by 
$\CK[\tau^{\pm1}]=\CK[\tau_i^{\pm1}(i\in I)]$ 
the ring of Laurent polynomials in the indeterminates 
$\tau_i$  ($i\in I$) with coefficients in $\CK$. 

\WSCI{\par\medskip}
\begin{thm} \label{thm-B}
Extend the automorphisms $s_i$ $(i\in I)$ of $\CK$ to those 
of $\CK[\tau^{\pm 1}]$ by the following action 
on the $\tau$-functions$\,:$ 
\begin{equation}
s_i(\tau_j)=\tau_j\quad(i\ne j),\qquad
s_i(\tau_i)=\varphi_i\, \tau_i \prod_{k\in I} \tau_k^{-a_{ki}}. 
\end{equation}
Then these $s_i$ $(i\in I)$ again give 
a realization of the Weyl group $W(A)$ as a group of 
automorphisms of $\CK[\tau^{\pm1}]$. 
\end{thm}
(For a more intrinsic formulation of this theorem, see 
Remark at the end of Section 4.)
\par\medskip
The action of $s_i$ on $\tau$-functions, defined above, 
is a multiplicative analogue of that of $r_i$ on the 
{\em fundamental weights} $\Lambda_j$ (modulo null roots),
except for the factor $\varphi_i$. 
Let $L=\bigoplus_{i\in I }\BZ\Lambda_i$ be the 
free $\BZ$-submodule of $\Hom_{\BZ\!}(Q^\vee,\BZ)$, 
generated by the dual basis $(\Lambda_i)_{i\in I}$ 
of $(h_i)_{i\in I}$. 
The Weyl group $W(A)=\br{\,r_i (i\in I) }$ act on $L$
so that 
\begin{equation}
r_i(\Lambda_j)=\Lambda_j \quad(i\ne j),\qquad
r_i(\Lambda_i)=\Lambda_i - \sum_{k\in I}\Lambda_k a_{k,i}.
\end{equation}
Note that 
$\br{w.h,w.\Lambda}=\br{h,\Lambda}$ 
($h\in Q^\vee, \Lambda\in L$) for any $w\in W(A)$ and that 
the natural homomorphism $Q\to L$ is a $W(A)$-homomorphism. 

In what follows we denote by $W=\br{\,s_i\,(i\in I)}$ the 
Weyl group $W(A)$ on the generators $s_i$ ($i\in I$). 
We introduce the notation of formal exponentials by setting 
\begin{equation}
\tau^\Lambda=\prod_{i\in I} \, \tau_i^\br{h_i,\Lambda}\qquad
(\Lambda\in L), 
\end{equation}
so that $\tau_i=\tau^{\Lambda_i}$ for $i\in I$. 
Then the ring $\CK[\tau^{\pm1}]$ of Laurent polynomials 
is isomorphic to the group ring $\CK[L]$.  
The action of $s_i$ on the $\tau$-functions, 
as defined in Theorem \ref{thm-B}, is rewritten in the form
\begin{equation}
s_i(\tau^{\Lambda})=\varphi_i^{\br{h_i,\Lambda}}\,\tau^{r_i.\Lambda}
\qquad(\Lambda\in L). 
\end{equation}
This implies that, for any $w\in W$ and $\Lambda\in L$, there 
exists a unique element $\phi_{w}(\Lambda) \in \CK$ such that
\begin{equation}
w(\tau^\Lambda)=\phi_{w}(\Lambda)\,\tau^{w.\Lambda}.
\end{equation}
Note that these $\phi_{w}(\Lambda)$ are determined by the
following cocycle condition:
\begin{equation}\begin{array}{l}
\smallskip
\phi_{w}(\Lambda+\Lambda')=\phi_{w}(\Lambda)\phi_w(\Lambda')
\quad(\Lambda,\Lambda'\in L),\cr
\smallskip
\phi_1(\Lambda)=1,\quad 
\phi_{s_i}(\Lambda)=\varphi_i^{\br{h_i,\Lambda}}\ \ (i\in I),\cr
\phi_{ww'}(\Lambda)= w(\phi_{w'}(\Lambda))\,\phi_w(w'.\Lambda)
\quad(w,w'\in W). 
\end{array}\end{equation}
This family 
$\phi=(\phi_{w}(\Lambda))_{w\in W,\Lambda\in L}$
of elements of $\CK$ will be called the {\em $\tau$-cocycle}.  
In fact, the correspondence 
\begin{equation}
\phi : W \to \Hom(L,\CK^{\times}) : 
w \mapsto (\Lambda \mapsto \phi_w(\Lambda))
\end{equation}
defines a 1-cocycle of the Weyl group $W(A)$ with coefficients 
in the $W(A)$-bimodule $\Hom(L,\CK^{\times})$.

The $\tau$-cocycle determines completely the action of the 
Weyl group on $\CK$.  
Note that, by the definition, the generators $\varphi_j$ for the 
Poisson algebra $\CA_0$ are expressed as 
\begin{equation}
\varphi_j=\frac{\tau_j \, s_j(\tau_j)}
{\prod_{i\ne j}\,\tau_i^{-a_{ij}}}\qquad(j\in I),
\end{equation}
multiplicatively in terms of $\tau$-functions.
This implies that 
\begin{equation}
w(\varphi_j)
=\frac{\phi_{w}(\Lambda_j) \, \phi_{ws_j}(\Lambda_j)}
{\prod_{i\ne j}\,\phi_{w}(\Lambda_i)^{-a_{ij}}}
\end{equation}
for any $j\in I,\ \ w\in W$.

By the definition of the $\tau$-cocycle, each 
$\phi_w(\Lambda_j)$ is {\em a priori} an element of the field of fractions 
$\CK=Q(\CA)$ of $\CA=\CA_0[\lambda]$. 
It turns out, however, that
the $\tau$-cocycle has a remarkable regularity 
in the following sense.  

\WSCI{\par\medskip}
\begin{thm}\label{thm-C}
Suppose that $A=(a_{ij})_{ij\in I}$ is a symmetrizable GCM. 
Then, 
for any $w\in W$ and $j\in I$, one has $\phi_{w}(\Lambda_j)\in\CA$,
namely,  $\phi_{w}(\Lambda_j)$ is a polynomial in $(\lambda_i)_{i\in I}$ 
with coefficients in $\CA_0$. 
\end{thm}
\WSCI{\par\medskip}
\noindent
Theorem implies that, if $\Lambda\in L$ is {\em dominant}, i.e., 
$\br{h_i,\Lambda}\ge 0$ for any 
$i\in I$, then one has $\phi_{w}(\Lambda)\in\CA=\CA_0[\lambda]$
for any $w\in W$. 
When $\CA_0$ has a {\em $\BZ$-form} in an appropriate sense, 
one can also show that $\phi_{w}(\Lambda)$ are 
defined over $\BZ$. 
(See Remark at the end of Section 5.)

These $\phi_{w}(\Lambda)$ can be regarded as a generalization 
of the so-called {\em Umemura polynomials} for generic solutions 
of the Painlev\'e equations\WSCI{.\cite{NOOU} }\PREP{ (\cite{NOOU}). }
It would be an interesting problem to investigate combinatorial 
properties of $\phi_w(\Lambda)$.

\rem{
Suppose that the generalized Cartan matrix $A$ is symmetrizable, 
and take nonzero rational numbers $\epsilon_i$ ($i\in I$) such 
that $a_{ij}\epsilon_j=a_{ji}\epsilon_{i}$ ($i,j\in I$). 
In this case, 
the realization of the Weyl group discussed in our previous paper \cite{NY} 
can be recovered essentially from the construction 
of this section by a special choice of the nilpotent Poisson 
algebra $\CA_0$.
For $\CA_0$, take the Poisson algebra,
{\em truncated at height 2}, defined by 
the generators $\varphi_i$ ($i\in I$) with the 
Serre relations 
\begin{equation}
\pad(\varphi_i)^{-a_{ij}+1}(\varphi_j)=0 \qquad(i\ne j),
\end{equation}
and 
\begin{equation}
\pbr{\varphi_i,\pbr{\varphi_j,\varphi_k}}=0 \qquad(i,j,k\in I). 
\end{equation}
To be consistent with the previous notation \cite{NY}, 
set $\alpha_i=\lambda_i/\epsilon_i$ and 
$u_{ij}=\epsilon_{i}\pbr{\varphi_i,\varphi_j}$. 
Since $\pbr{u_{ij},\varphi}=0$ for any $\varphi\in\CA_0$, 
one can treat $u_{ij}$ as constants.
}

\section{Examples of rank 2}
In this section we give examples of our realization of 
the Weyl groups for generalized Cartan matrices of rank 2.  
In the following examples, we use the notation
\begin{equation}
\alpha=\alpha_1,\quad \beta=\alpha_2,\quad
a=\lambda_1,\quad b=\lambda_2.
\end{equation}

\subsection{Case of $2 A_1$}  We set
$x=\varphi_1, y=\varphi_2$ and $\pbr{x,y}=0$. 
In this case, our realization of the Weyl group is 
trivial on the variables $x$, $y$.
\begin{equation}
\begin{picture}(80,28)(-20,10)
\put(0,0){\vector(1,0){30}}
\put(0,0){\vector(-1,0){30}}
\put(30,2){\small $\alpha:x$}
\put(0,0){\vector(0,1){30}}
\put(5,28){\small $\beta:y$}
\end{picture}
\quad
A=\left[\begin{array}{cc}
2 & 0 \\ 0 &2
\end{array}\right]
\qquad
\begin{array}{c||cc|cc}
 & a  & b& x & y \cr
\hline
s_1 &-a &b &x & y\cr
s_2 &a &-b &x & y 
\end{array}
\end{equation}

\subsection{Case of $A_2$} 
We set 
$x=\varphi_1$, $y=\varphi_2$, and 
$z=\pbr{x,y}=\pbr{\varphi_1,\varphi_2}. $
\begin{equation}
\begin{picture}(80,30)(-20,10)
\put(0,0){\vector(1,0){30}}
\put(0,0){\vector(-1,0){30}}
\put(30,2){\small $\alpha:x$}
\put(0,0){\vector(-2,3){15}}
\put(-20,25){\small $\beta:y$}
\put(0,0){\vector(2,3){15}}
\put(15,25){\small $\alpha+\beta:z$}
\end{picture}
\quad
A=\left[\begin{array}{cc}
2 & -1 \\ -1 &2 
\end{array}\right]
\qquad
\begin{array}{c|ccc}
\pbr{\,,\,} & x & y & z\\
\hline
x & 0 & z & 0\\
y & -z & 0 & 0\\
z & 0 & 0 & 0 
\end{array}
\end{equation}
\par\medskip\noindent 
We take the Poisson algebras $\CA_0=\BC[x,y,z]$ 
(or any quotient Poisson algebra of $\BC[x,y,z]$ without zerodivisors),
and set $\CA=\BC[a,b,x,y,z]$.  
By the table of the Poisson bracket indicated above, 
the action of $s_1$, $s_2$ is determined as follows:
\begin{equation}
\begin{array}{c||cc|ccc}
 & a  & b& x & y & z\\
\hline
s_1 &-a &b+a &x & y+\frac{a z}{x} & z \\
s_2 &a+b &-b &x-\frac{b z}{y} & y & z
\end{array}
\end{equation}
The automorphisms $s_1, s_2$ defined by this table give a 
realization of the Weyl group $ W(A_2)=\br{s_1,s_2} \iso \frS_3$ 
as a group of canonical transformations on the 
field of rational functions $\CK=Q(\CA)=\BC(a,b,x,y,z)$. 

The $\tau$-cocycle $\phi=(\phi_{w}(\Lambda))_{w,\Lambda}$ 
for this realization is given as follows. 
\begin{equation}
\begin{array}{c|cccccc}
\Lambda\backslash w
 & \ 1\ & s_1 & \ s_2 \ & \ s_2s_1\ &\ s_1s_2\ &\ s_1s_2s_1=s_2s_1s_2\ 
\\[2pt]\hline
\ \Lambda_1\  & 1 & x & 1 & xy-bz & x & xy-bz
\cr
\Lambda_2 & 1 & 1 & y &  y  & xy+az & xy+az
\end{array}
\end{equation}
We remark that
the value $\phi_{w_0}(\Lambda_1)$
of the $\tau$-cocyle at the longest element $w_0$ can be
determined in two ways as
\begin{equation}
\phi_{s_2s_1s_2}(\Lambda_1)=\varphi_2s_2(\varphi_1),\quad
\phi_{s_1s_2s_1}(\Lambda_1)=s_1(\varphi_2s_2(\varphi_1)),
\end{equation}
by using the cocycle property for the two reduced 
decompositions of 
$w_0$. 
One can verify that $\varphi_2s_2(\varphi_1)
=y(x-\dfrac{bz}{y})=xy-b z$ is invariant with respect to $s_1$,
which guarantees the equality of the two expressions above. 

\subsection{Case of $B_2$} 
We set $x=\varphi_1,\  y=\varphi_2$ and 
\begin{equation}
 z=\pbr{x,y}=\pbr{\varphi_1,\varphi_2},\quad
w=\frac{1}{2}\pbr{x,z}=\frac{1}{2}\pbr{\varphi_1,\pbr{\varphi_1,\varphi_2}}.
\end{equation}
\begin{equation}
\begin{picture}(100,35)(-20,10)
\put(0,0){\vector(-1,1){24}}
\put(0,0){\vector(0,1){24}}
\put(0,0){\vector(-1,0){24}}
\put(0,0){\vector(1,0){24}}
\put(0,0){\vector(1,1){24}}
\put(30,0){\small $\alpha:x$}
\put(-44,28){\small $\beta:y$}
\put(-15,28){\small $\alpha+\beta:z$}
\put(28,28){\small $2\alpha+\beta:w$}
\end{picture}
\qquad
A=\left[\begin{array}{cc}
2 & -2 \\ -1 &2 
\end{array}\right]
\end{equation}
\par\medskip\noindent
The Poisson bracket and the action of $s_1$, $s_2$ are 
defined by the following tables:
\begin{equation}
\begin{array}{c|cccc}
\pbr{\,,\,}  & x & y & z & w\\
\hline
x & 0 & z & 2w & 0\\
y & -z & 0 & 0 & 0\\
z & -2w & 0 & 0 &0 \\
w& 0 & 0 & 0 & 0 
\end{array}
\end{equation}
$$
\begin{array}{c||cc|cccc}
 &a & b & x & y & z &w \\
\hline
s_1&-a &a+b &x & y+\frac{a z}{x}+\frac{a^2 w}{x^2} & 
z+\frac{2aw}{x} & w\\
s_2 &a+2 b &-b &x-\frac{b z}{y} 
& y & z & w
\end{array}
$$
From the nilpotent Poisson algebra $\CA_0=\BC[x,y,z,w]$, 
we obtain the realization of the Weyl group $W(B_2)$ as 
a group of canonical transformations on the field of 
rational functions $\CK=\BC(a,b,x,y,z,w)$. 

The $\tau$-cocycle $\phi=(\phi_{w}(\Lambda))_{w,\Lambda}$ 
for this realization is given as follows. 
\begin{equation}
\begin{array}{c|ccccc}
 & \ 1\ & s_1 & \ s_2 \ & \ s_2s_1\ &\ s_1s_2\ 
\\[2pt]\hline
\ \Lambda_1\  & 1 & x & 1 & xy-bz & x 
\cr
\Lambda_2 & 1 & 1 & y &  y  & \ x^2 y+a\,x z+a^2\,w \ 
\end{array}
\end{equation}
$$
\begin{array}{c|ccc}
 & \ s_1 s_2 s_1 \ & \ s_2 s_1 s_2 
\\[2pt]\hline
\ \Lambda_1\  & \ x^2 y -b\,x z -a(a+2b)\,w \ 
& xy-bz 
\cr
\Lambda_2 &  x^2 y+a\,x z+a^2\,w & x^2 y^2 +a\,x y z
-(a+b) b\,z^2+(a+2b)^2 y w
\end{array}
$$ 
$$
\begin{array}{c|c}
 & \ s_2 s_1 s_2 s_1=s_1s_2 s_1 s_2 
\\[2pt]\hline
\ \Lambda_1\  & \ x^2 y -b\,x z -a(a+2b)\,w \ 
\cr
\Lambda_2 &  
\ x^2 y^2 +a\,x y z-(a+b) b\,z^2+(a+2b)^2 y w\ 
\end{array}
$$
We remark that 
\begin{equation}
\phi_{s_1s_2s_1}(\Lambda_1)=\varphi_1 s_1(\varphi_2\,s_2(\varphi_1))
=x^2 y -b\,x z -a(a+2b)\,w
\end{equation}
is $s_2$-invariant.  From this fact, it follows that the 
two expressions 
\begin{equation}
\phi_{s_1s_2s_1s_2}(\Lambda_1)=\phi_{s_1s_2s_1}(\Lambda_1), 
\quad
\phi_{s_2s_1s_2s_1}(\Lambda_1)=s_2(\phi_{s_1s_2s_1}(\Lambda_1))
\end{equation}
give the same value $\phi_{w_0}(\Lambda_1)$ for the longest element 
$w_0$ of the Weyl group. 
Similarly, $\phi_{w_0}(\Lambda_2)$ is determined consistently 
from the $s_1$-invariance of 
\begin{equation}
\begin{array}{l}
\smallskip
\phi_{s_2s_1s_2}(\Lambda_2)=
\varphi_2 s_2(\varphi_1^2\,s_1(\varphi_2))\cr
\phantom{\phi_{s_2s_1s_2}(\Lambda_2)}
= x^2 y^2 +a\,x y z-(a+b) b\,z^2+(a+2b)^2 y w. 
\end{array}
\end{equation}

\subsection{Case of $G_2$} 
We set $u=\varphi_1,\  v=\varphi_2$ and 
\begin{equation}
 w=\pbr{u,v},\quad x=\frac{1}{2}\pbr{u,w}, \quad y=\frac{1}{3}\pbr{u,x}, 
\quad z=\pbr{v,y}.
\end{equation}
\begin{equation}
\begin{picture}(140,60)(-60,10)
\put(0,0){\vector(2,1){50}}
\put(0,0){\vector(2,3){18}}
\put(0,0){\vector(-2,3){18}}
\put(0,0){\vector(1,0){36}}
\put(0,0){\vector(-1,0){36}}
\put(0,0){\vector(0,1){50}}
\put(0,0){\vector(-2,1){50}}
\put(40,0){\small $\alpha:u$}
\put(55,20){\small $3\alpha+\beta:y$}
\put(12,32){\small $2\alpha+\beta:x$}
\put(-45,32){\small $\alpha+\beta:w$}
\put(-75,20){\small $\beta:v$}
\put(-15,55){\small $3\alpha+2\beta:z$}
\end{picture}
\qquad\quad
A=\left[\begin{array}{cc}
2 & -3 \\ -1 &2 
\end{array}\right]
\end{equation}
\par\medskip\noindent
The Poisson bracket and the action of $s_1$, $s_2$ are determined 
by the following tables. 
\begin{equation}
\begin{array}{c|cccccc}
\pbr{\,,\,}  & u & v & w& x & y & z\\
\hline
u & 0 & w & 2x & 3y & 0 & 0 \\
v & -w & 0 & 0 & 0 & z & 0\\
w & -2x & 0 & 0 &-3z & 0 & 0 \\
x & -3y & 0 & 3z & 0 & 0 & 0\\ 
y & 0 & -z & 0 & 0 & 0 & 0\\
z & 0 & 0 & 0 & 0  & 0 & 0
\end{array}
\qquad
\begin{array}{c||cc}
 &a & b   \\
\hline
s_1 &-a &a+b \\
s_2 &a+3 b &-b 
\end{array}
\end{equation}
$$
\begin{array}{c||cccccc}
 &u & v & w & x & y & z  \\
\hline
s_1 & u &
v+\frac{aw}{u}+ \frac{a^2x}{u^2}+ \frac{a^3y}{u^3} &
w+\frac{2ax}{u}+\frac{3a^2y}{u^2}& 
x +\frac{3ay}{u} &
y & z\\
s_2 & u-\frac{bw}{v} & v & w & x & y+\frac{b z}{v} & z
\end{array}
$$
The nilpotent Poisson algebra $\CA_0=\BC[u,v,w,x,y,z]$
defines a group of canonical transformations, isomorphic 
to the Weyl group $W(G_2)$, on the field of rational functions 
$\CK=\BC(a,b,u,v,w,x,y,z)$. 

The $\tau$-cocycle in this case is given as follows. 
$$
\begin{array}{l}
\smallskip
\phi_1(\Lambda_1)=\phi_{s_2}(\Lambda_1)=1,\qquad
\phi_{s_1}(\Lambda_1)=\phi_{s_1s_2}(\Lambda_1)
=u\cr\smallskip
\phi_{s_2s_1}(\Lambda_1)=\phi_{s_2 s_1s_2}(\Lambda_1)
=uv-bw,\cr\smallskip
\phi_{s_1s_2s_1}(\Lambda_1)=\phi_{s_1s_2 s_1s_2}(\Lambda_1)
=u^3 v-b\,u^2 w-a(a+2b)\,ux-a^2(2a+3b)y\cr\smallskip
\phi_{s_2s_1s_2s_1}(\Lambda_1)=\phi_{s_2s_1s_2 s_1s_2}(\Lambda_1)
\cr\smallskip
\quad
=u^3v^2-2b\,u^2 v w 
-(a+b)(a+3b)\,uvx
+b^2\,u w^2
\cr\smallskip
\qquad
-(a+3b)^2(2a+3b)\,vy
+b(a+b)(a+3b)\,wx
-b(a+3b)^2(2a+3b)\,z
\cr\smallskip
\phi_{s_1s_2s_1s_2s_1}(\Lambda_1)=\phi_{s_1s_2s_1s_2 s_1s_2}(\Lambda_1)
\cr\smallskip
\quad=
u^4 v^2-2 b u^3 v w -(2 a^2+6 a b+3 b^2)u^2 vx+b^2 u^2 w^2
-(2a+3b)^3 uvy
\cr\smallskip\qquad
+b(2a^2+6ab+3b^2)uwx
-(a+b)(a+3b)(2a+3b)^2uz
\cr\smallskip\qquad
-a(a+2b)(2a+3b)^2wy
+a(a+b)(a+2b)(a+3b)x^2
\PREP{
\end{array}
$$
\vspace{-20pt}
$$
\begin{array}{l}
}
\\[6pt]
\smallskip
\phi_{1}(\Lambda_2)=\phi_{s_1}(\Lambda_2)=1,
\quad
\phi_{s_2}(\Lambda_2)=\phi_{s_2s_1}(\Lambda_2)=v
\cr\smallskip
\phi_{s_1s_2}(\Lambda_2)=\phi_{s_1s_2s_1}(\Lambda_2)
=u^3 v  +a\,u^2w +a^2\,u x+a^3\,y
\cr\smallskip
\phi_{s_2s_1s_2}(\Lambda_2)=\phi_{s_2s_1s_2s_1}(\Lambda_2)
\cr\smallskip
\quad
=u^3 v^3
+a\,u^2v^2w
+(a+3b)^2\,uv^2 x
-b(2a+3b)\,uvw^2
\cr\smallskip
\qquad
+(a+3b)^3\,v^2y
-b(a+3b)^2\,vwx
+b(a+3b)^3 vz
+b^2(a+2b)\,w^3
\WSCI{\cr\smallskip}
\PREP{
\end{array}$$\quad
$$
\begin{array}{l}
}
\phi_{s_1s_2s_1s_2}(\Lambda_2)=\phi_{s_1s_2s_1s_2s_1}(\Lambda_2)
\cr\smallskip
\quad
=
u^6v^3 + 2 a\,u^5v^2w 
+(5a^2+12a+9b^2)\,u^4v^2x
-b(4a+3b)\,u^4vw^2
\cr\smallskip
\qquad
+(2a+3b)^2(5a+3b)\,u^3 v^2 y
-b(8a^2+15ab+9b^2)\,u^3 v w x
\cr\smallskip
\qquad
+(a+b)(2a+3b)^3\,u^3 v z
+2b^2(a+b)\,u^3 w^3
+a(2a+3b)^2(5a+3b)\,u^2 v w y
\cr\smallskip\qquad
-a(5a^3+24a^2b+36ab^2+18b^3)\,u^2 v x^2
+3ab^2(a+b)\,u^2 w^2x
\cr\smallskip\qquad
+a(a+b)(2a+3b)^3\,u^2 w z
-a^2(a+3b)(2a+3b)^2\,u v x y
\cr\smallskip\qquad
+2a^2(a+b)(2a+3b)^2\,u w^2 y
-a^2(a+b)(2a^2+6ab+3b^2)\,u w x^2
\cr\smallskip\qquad
+a^2(a+b)(2a+3b)^2\,x z
-a^3(2a+3b)^3\, v y^2
+a^3(a+b)(2a+3b)^2\,w x y
\cr\smallskip\qquad
-a^3(a+b)^2(a+2b)\,x^3
+a^3(a+b)(2a+3b)^3\,yz
\end{array}
$$
$$
\begin{array}{l}
\smallskip
\phi_{s_2s_1s_2s_1s_2}(\Lambda_2)=\phi_{s_2s_1s_2s_1s_2s_1}(\Lambda_2)
\cr\smallskip
\quad
=
u^6 v^4+2au^5v^3w+(5a^2+18ab+18b^2)u^4v^3x-6b(a+b)u^4v^2w^2
\cr\smallskip\quad
+(2a+3b)^2(5a+12b)u^3v^3y-3b(4a^2+13ab+12b^2)u^3v^2wx
\cr\smallskip \quad
+2(a+3b)^2(2a+3b)^2u^3v^2z
+2b^2(3a+4b)u^3vw^3
\cr\smallskip\quad
+a(2a+3b)^2(5a+12b)u^2v^2wy 
-(a+3b)(5a^3+21a^2b+27ab^2+9b^3)u^2v^2x^2
\cr\smallskip\quad
+3b^2(3a^2+8ab+6b^2)u^2vw^2x
+2a(a+3b)^2(2a+3b)^2u^2vwz
\cr\smallskip\quad
-b^3(2a+3b)u^2w^4
-a(a+3b)^2(2a+3b)^2uv^2xy
\cr\smallskip\quad
+a(2a+3b)^2(2a^2+6ab+3b^2)uvw^2y
\cr\smallskip\quad
-a(a+3b)(2a^2+6a^2b+3ab^2-3b^3)uvwx^2
\cr\smallskip\quad
+2(a+3b)^2(2a+3b)^2(a^2+3ab+3b^2)uvxz
-ab^3(2a+3b)uw^3x
\cr\smallskip\quad
-2b(a+b)(a+3b)^2(2a+3b)^2uw^2z -(a+3b)^3(2a+3b)^3v^2y^2
\cr\smallskip\quad
+ (a+3b)^2(2a+3b)^2(a^2+6ab+6b^2)vwxy 
-(a+b)(a+2b)^2(a+3b)^3vx^3 
\cr\smallskip\quad
+a(a+3b)^3(2a+3b)^3vyz
-b(a+2b)^2(2a+3b)^3 w^3y 
\cr\smallskip\quad
+b(a+b)(a+3b)(2a+3b)(a^2+3ab+3b^2)w^2x^2
\cr\smallskip\quad
-ab(a+b)(a+3b)^2(2a+3b)^2wxz
+b(a+b)(a+3b)^3(2a+3b)^3z^2
\end{array}
$$
We remark that 
\begin{equation}
\begin{array}{l}\smallskip
\phi_{s_1s_2s_1s_2s_1}(\Lambda_1)=
\varphi_1s_1\varphi_2s_2\varphi_1^2s_1\varphi_2s_2(\varphi_1),
\quad\mbox{and}
\cr
\phi_{s_2s_1s_2s_1s_2}(\Lambda_2)=
\varphi_2 s_2 \varphi_1^3 s_1 \varphi_2^2 s_2 \varphi_1^3 s_1(\varphi_2)
\end{array}
\end{equation}
are invariant with respect to 
$s_2$ and $s_1$, respectively. 

\PREP{\newpage}

\section{Realization of  Weyl groups}

In this section, we give a proof of Theorems \ref{thm-A}
and \ref{thm-B}. 
In fact we prove these two theorems as a consequence of 
their {\em formal} version (Theorem \ref{thm-E} below).  

\subsection{Formal version}
Fixing a GCM $A=(a_{ij})_{i,j\in I}$, 
let $\CR$ be a Poisson which contains 
a set of invertible elements $(\varphi_i)_{i\in I}$ such that 
\begin{equation}
\pad(\varphi_i)^{-a_{ij}+1}(\varphi_j)=0\quad(i\ne j). 
\end{equation}
We denote by $\CR[[\lambda]]$ the ring of formal power series 
in the $\lambda$-variables $(\lambda_i)_{i\in I}$ with 
coefficients in $\CR$.  
The Poisson bracket on $\CR$ extends naturally to $\CR[[\lambda]]$ 
by the trivial action on the $\lambda$-variables. 
We also use the same notation 
$r_i :\CR[[\lambda]] \to \CR[[\lambda]]$
for the formal completion of the 
$\CR$-linear automorphism $r_i: \CR[\lambda]\to\CR[\lambda]$
such that $r_i(\lambda_j)=\lambda_j-a_{ji}\lambda_i$ ($j\in I$), 
For each $i\in I$, we define the linear mapping 
$s_i : \CR[[\lambda]]\to \CR[[\lambda]]$ by 
\begin{equation}
s_i=\exp\left(\frac{\lambda_i}{\varphi_i}\pad(\varphi_i)\right)\circ r_i
\quad(i\in I). 
\end{equation}
In the following we set
\begin{equation}
X_i=\frac{1}{\varphi_i}\pad(\varphi_i)
\qquad(i\in I),
\end{equation}
so that $R_i(t)=\exp(t X_i)$ and $s_i=R_i(\lambda_i)\circ r_i$. 

\WSCI{\par\medskip}
\begin{lem} 
For each $i\in I$, $s_i$ defines a canonical transformation of 
$\CR[[\lambda]]$ $:$
\begin{equation}
s_i(\varphi\psi)=s_i(\varphi)s_i(\psi), \quad
s_i(\pbr{\varphi,\psi})=\pbr{s_i(\varphi), s_i(\psi)}, 
\end{equation}
for any $\varphi, \psi\in \CR[[\lambda]]$. 
\end{lem}
\WSCI{\par\medskip}
\noindent
This lemma is a consequence of the fact that the derivations 
$X_i$ have the property
\begin{equation}
X_i(\pbr{f,g})=\pbr{X_i(f),g}+\pbr{f,X_i(g)}\qquad(f,g\in \CR).
\end{equation}

We denote by $\hW=\br{ s_i (i\in I)}$ the free group 
generated by the symbols $s_i$ ($i\in I$), and 
by $\rho: \hW\to W$ the homomorphism defined by 
$\rho(s_i)=r_i$ ($i\in I$). 
At this stage, we have a group homomorphism 
$\hW \to \Aut(\CR[[\lambda]])$, and this action, 
restricted to the $\lambda$-variables,
is factored through $\rho: \hW\to W$. 

We also consider the ring of Laurent polynomials 
$\CR[[\lambda]][\tau^{\pm1}]$
in the $\tau$-variables, and extend the action of $s_i$ 
to $\CR[[\lambda]][\tau^{\pm1}]$ by
\begin{equation}
s_i(\tau_j)=\tau_j  \ \ (i\ne j),\quad
s_i(\tau_i)=\varphi_i\,\tau_i \prod_{k\in I}\tau_k^{-a_{ki}}. 
\end{equation}
In this setting, we have a group homomorphism 
$\hW \to \Aut(\CR[[\lambda]][\tau^{\pm1}])$. 
Our goal is to prove 

\WSCI{\par\medskip}
\begin{thm}\label{thm-E}
The automorphisms $s_i $ $(i\in I)$ of 
$\CR[[\lambda]][\tau^{\pm1}]$, 
defined as above, satisfy the fundamental relations 
for the Weyl group $W(A)$ $:$
\begin{equation}
s_i^2=1, \quad (s_i\,s_j)^{m_{ij}}=1\quad(i\ne j).
\end{equation} 
\end{thm}
\WSCI{\par\medskip}

We show how Theorem \ref{thm-E} implies Theorems 
\ref{thm-A} and \ref{thm-B}.  
With the notation of Section 1, we take the localization 
$\CR=\CA_0[\,\varphi_j^{-1} (j\in I)]$. 
In order to establish Theorem \ref{thm-A}, we 
have only to check the validity of equalities 
$s_i^2(\psi)=1$ and $(s_{i}s_{j})^{m_{ij}}(\psi)=1$ 
for any $\psi\in\CA_0$.  
By the definition of $s_j$ ($j\in I$), the elements 
$s_i^2(\psi)$, $(s_{i}s_{j})^{m_{ij}}(\psi)$ of 
$\CK=Q(\CA_0)(\lambda)=Q(\CR)(\lambda)$  
are regular at $\lambda=0$.  
Hence the equalities 
$s_i^2(\psi)=1$ and $(s_{i}s_{j})^{m_{ij}}(\psi)=1$ 
follow from those in the ring of formal 
power series $\CR[[\lambda]]$ which are guaranteed 
by Theorem \ref{thm-E}. 
Also, Theorem \ref{thm-B}, concerning the 
action of $s_i$ on $\CK[\tau^{\pm1}]$, 
follows similarly from 
the equalities $s_i^2(\tau_j)=\tau_j$ and 
$(s_{i}s_{j})^{m_{ij}}(\tau_j)=\tau_j$ in 
$\CR[[\lambda]][\tau^{\pm 1}]$. 

\subsection{$\tau$-Cocycle and braid relations on $\tau$-functions}

We first check the equality $s_i^2=1$ on $\CR[[\lambda]]$ 
for each $i\in I$.   
In fact, with the $R$-operator notation $R_i(t)=\exp(t X_i)$, one has
\begin{equation}
s_i^2=R_i(\lambda_i) r_i R_i(\lambda_i) r_i
=R_i(\lambda_i) R_i(-\lambda_i) r_i^2=1,
\end{equation}
since $R_i(-t)=\exp(-t X_i)=R_i(t)^{-1}$. 
It is also easy to see $s_i^2(\tau_j)=\tau_j$ for all $j\in I$. 
\par\medskip
Instead of the Coxeter relations $(s_i s_j)^{m_{ij}}=1$, 
we will verify the following braid relations:
\begin{equation}\label{BR}
\begin{array}{clcc}
(0) & (a_{ij},a_{ji})=(0,0) &\Longrightarrow& s_is_j=s_js_i,\cr
(1) & (a_{ij},a_{ji})=(-1,-1) &\Longrightarrow& s_is_js_i=s_js_is_j,\cr
(2) & (a_{ij},a_{ji})=(-2,-1) &\Longrightarrow& s_is_js_is_j=s_js_is_js_i,\cr
(3) & (a_{ij},a_{ji})=(-3,-1) &\Longrightarrow& 
s_is_js_is_js_is_j=s_js_is_js_is_js_i,
\end{array}
\end{equation}
which corresponds to the root systems of rank 2 
of type $2A_1$, $A_2$, $B_2$ and $G_2$, respectively. 
Note that the validity of each of these braid relations on $\CR[[\lambda]]$
is equivalent to the corresponding Yang-Baxter equation
in (\ref{YBE}) for $R_j(t)$.  
Since the commutativity $s_is_j=s_js_i$ for the case $a_{ij}=a_{ji}=0$ 
is immediate, we will mainly consider the other three cases. 

\WSCI{\par\medskip}
\begin{lem} \label{lem-A}
The following identities hold in $\CR[[\lambda]]$. 
\smallskip
\newline
$(1)$  If $(a_{ij},a_{ji})=(-1,-1)$, 
\begin{equation}
\varphi_j s_j(\varphi_i)=s_i \varphi_j s_j(\varphi_i),
\quad
\varphi_i s_i(\varphi_j)=s_j \varphi_i s_i(\varphi_j).
\end{equation}
$(2)$  If $(a_{ij},a_{ji})=(-2,-1)$, 
\begin{equation}
\varphi_i s_i \varphi_j s_j(\varphi_i)=
s_j\varphi_i s_i\varphi_js_j(\varphi_i),\quad
\varphi_j s_j \varphi_i^2 s_i(\varphi_j)=
s_i\varphi_js_j\varphi_i^2 s_i(\varphi_j).
\end{equation}
$(3)$ If $(a_{ij},a_{ji})=(-3,-1)$, 
\begin{equation}
\begin{array}{c}
\smallskip
\varphi_is_i \varphi_j s_j\varphi_i^2
s_i\varphi_j s_j(\varphi_i)
=s_j \varphi_is_i\varphi_js_j\varphi_i^2
s_i\varphi_js_j(\varphi_i),\\
\varphi_js_j\varphi_i^3s_i\varphi_j^2 
s_j\varphi_i^3s_i(\varphi_j)
=s_i\varphi_js_j\varphi_i^3s_i\varphi_j^2
s_j\varphi_i^3s_i(\varphi_j).
\end{array}
\end{equation}
\end{lem}
\WSCI{\par\medskip}
\noindent
In this lemma, we regard $\varphi_i$ and $\varphi_j$ as 
multiplication operators. 
The two identities in the case of $(a_{ij},a_{ji})=(-1,-1)$,
for example, should be read as 
\begin{equation}
\varphi_j s_j(\varphi_i)=s_i(\varphi_j) s_i s_j(\varphi_i),
\quad
\varphi_i s_i(\varphi_j)=s_j(\varphi_i) s_j s_i(\varphi_j).
\end{equation}
We will not write down the corresponding formulas for the other 
two cases.  
Note that, in each case, the exponents of $\varphi_i, \varphi_j$ in 
the two formulas correspond to the coefficients of $a$ and $b$ 
in the Yang-Baxter equation (\ref{YBE}), respectively. 
The six equalities in 
Lemma \ref{lem-A} can be verified by 
direct computation,   
by using the automorphisms $s_1$, $s_2$ as in the examples 
of the previous section. 
(In fact, we made use of a computer algebra system to check 
Lemma \ref{lem-A} for $B_2$ and $G_2$.)
We will explain now how Lemma \ref{lem-A} imply the Yang-Baxter 
equations (\ref{YBE}) and the braid relations (\ref{BR}). 

\par\medskip
In order to clarify the meaning of Lemma \ref{lem-A}, 
we use the terminology of {\em $\tau$-cocycle}. 
With the notation of formal exponentials $\tau^\Lambda$ 
($\Lambda\in L$), for each $w\in \hW$, one can define 
the elements 
$\phi_{w}(\Lambda)\in \CR[[\lambda]]$ ($\Lambda\in L$)
by the formula  
\begin{equation}
w(\tau^\Lambda) =\phi_w(\Lambda)\,\tau^{\rho(w).\Lambda}
\quad(w\in\hW, \Lambda\in L). 
\end{equation}
By the definition, these $\phi_w(\Lambda)$ have the following cocycle property. 

\WSCI{\par\medskip}
\begin{prop}
$(1)$ 
The elements $\phi_w(\Lambda)$ are determined uniquely by  
\begin{equation}
\begin{array}{l}
\smallskip
\phi_1(\Lambda)=1,\quad
\phi_{s_i}(\Lambda)=\phi_{s_i^{-1}}(\Lambda)
=\varphi_i^{\br{h_i,\Lambda}}\quad(i\in I),
\\
\phi_w(\Lambda+\Lambda')=\phi_w(\Lambda)\phi_w(\Lambda'),\quad\phi_{ww'}(\Lambda)=
\phi_{w}(\rho(w').\Lambda) w(\phi_{w'}(\Lambda)),
\end{array}
\end{equation}
for any $w,w'\in \hW$ and $\Lambda, \Lambda'\in L$.
\newline
$(2)$ For any $w=s_{j_1}s_{j_2}\cdots s_{j_p}\in \hW$ and $\Lambda\in L$, 
$\phi_w(\Lambda)$ is expressed as 
\begin{equation}\label{phi}
\begin{array}{l}
\smallskip
\phi_w(\Lambda)=\prod_{k=1}^p \, 
s_{j_1}\cdots s_{j_{k-1}}(\varphi_{j_k}^{\br{h_{j_k},r_{j_{k+1}}\cdots r_{j_{p}}\Lambda}})\\
=
\varphi_{j_1}^{\br{h_{j_1},r_{j_2}\cdots r_{j_p}\Lambda}}s_{j_1}
\varphi_{j_2}^{\br{h_{j_2},r_{j_3}\cdots r_{j_p}\Lambda}}\cdots\cr
\qquad \cdot \ s_{j_1}\cdots s_{j_{p-2}}
\varphi_{j_{p-1}}^{\br{h_{j_{p-1}},r_{j_p}\Lambda}}s_{j_1}\cdots
s_{j_{p-1}}(\varphi_{j_p}^{\br{h_{j_p},\Lambda}}).
\end{array}
\end{equation}
\end{prop}
\WSCI{\par\medskip}
\noindent
At this stage, we only know that 
$\phi$ is a 1-cocycle of $\hW$ with coefficients 
in the $\hW$-bimodule $\Hom(L,\CK^\times)$.

\par\medskip 
The formulas of Lemma \ref{lem-A} are understood as relations among 
$\phi_{w}(\Lambda)$ for certain $w\in \hW$ and $\Lambda\in L$. 
For example, let us take two indices $i,j\in I$ such that $(a_{ij},a_{ji})=(-1,-1)$.  
Then, as a special case of formula (\ref{phi}), we have 
\begin{equation}
\begin{array}{l}
\smallskip
\phi_{s_js_is_j}(\Lambda_i)=\varphi_j s_j(\varphi_i),\quad
\phi_{s_is_js_i}(\Lambda_i)=s_i \varphi_j s_j(\varphi_i),\\
\phi_{s_is_js_i}(\Lambda_j)=\varphi_i s_i(\varphi_j),\quad
\phi_{s_js_is_j}(\Lambda_j)=s_j \varphi_i s_i(\varphi_j).\
\end{array}
\end{equation}
Note also that $\phi_{s_is_js_i}(\Lambda_k)=1$ for $k\ne i,j$. 
Hence, Lemma \ref{lem-A} (1) implies
\begin{equation}
\phi_{s_is_js_i}(\Lambda_k)=\phi_{s_js_is_j}(\Lambda_k)\quad (k\in I).
\end{equation}
In the same way, Lemma \ref{lem-A} can be reformulated as follows. 

\WSCI{\par\medskip}
\begin{lem}\label{LemB}
For each $i,j\in I$, the following identities hold in $\CR[[\lambda]]$.
\smallskip\newline
$(1)$  If $(a_{ij},a_{ji})=(-1,-1)$, \quad
$\phi_{s_is_js_i}(\Lambda)=\phi_{s_js_is_j}(\Lambda)\quad (\Lambda\in L).$
\smallskip\newline
$(2)$  If $(a_{ij},a_{ji})=(-2,-1)$, \quad
$\phi_{s_is_js_is_j}(\Lambda)=\phi_{s_js_is_js_i}(\Lambda)\quad (\Lambda\in L).$
\smallskip\newline
$(3)$ If $(a_{ij},a_{ji})=(-3,-1)$, \quad
$\phi_{s_is_js_is_js_is_j}(\Lambda)=\phi_{s_js_is_js_is_js_i}(\Lambda)\quad (\Lambda\in L).$
\end{lem}
\WSCI{\par\medskip}
\noindent
Note that Lemma \ref{LemB} implies that the braid relations are 
valid on the $\tau$-functions.

\WSCI{\par\medskip}
\begin{cor}\label{cor-BRtau}
When  $a_{ij}a_{ji}=0,1,2,$ or $3$, one has
\begin{equation}\label{BRtau} 
s_is_js_i\cdots (\tau^\Lambda)=s_js_is_j\cdots(\tau^\Lambda)
\quad(\mbox{$m_{ij}$ factors of $s_k$'s on each side})
\end{equation}
for all $\Lambda\in L$.  
\end{cor}

\subsection{Braid relations on $\CR[[\lambda]]$}
We now proceed to the braid relations on the Poisson algebra $\CR[[\lambda]]$. 
For the argument below, it is convenient to extend the Poisson algebra 
$\CR$ to 
\begin{equation}
\widetilde{\CR}=\BC[\log \varphi_j \,(j\in I)]\otimes_{\BC} \CR,
\end{equation}
by adjoining the {\em formal logarithms} of the elements $\varphi_i$ ($i\in I$).
The Poisson bracket $\pbr{\,,\,}$ extends naturally to $\widetilde{\CR}$
and $\widetilde{\CR}[[\lambda]]$ by 
\begin{equation}
\begin{array}{c}
\pbr{\left(\log\varphi_j\right)^m,\psi}=m(\log\varphi_j)^{m-1}
\dfrac{1}{\varphi_j}\pbr{\varphi_j,\psi},\cr
\pbr{\left(\log\varphi_i\right)^m,\left(\log\varphi_j\right)^n}=
mn(\log\varphi_i)^{m-1}(\log\varphi_j)^{n-1}
\dfrac{1}
{\varphi_i\varphi_j}\pbr{\varphi_i,\varphi_j},
\end{array}
\end{equation}
for any $i,j\in I$, $\psi\in \widetilde{\CR}$ and $m,n=0,1,2,\ldots$. 
Consequently, one can define the action of 
$s_i$ on $\widetilde{\CR}[[\lambda]]$ by 
\begin{equation}
s_i=\exp(\lambda_i X_i)\circ r_i,\quad 
X_i=\pad(\log \varphi_i)=\frac{1}{\varphi_i}\pad(\varphi_i) \quad(i\in I).
\end{equation}

Note that each Yang-Baxter equation in (\ref{YBE}) takes the form 
\begin{equation}\label{EE1}
\begin{array}{l}
\smallskip
\exp(\mu_1 X_i) \exp(\mu_2 X_j)\cdots\exp(\mu_m X_\ell)\\
\quad=\ \ \exp(\mu_m X_j) \exp(\mu_{m-1} X_i)\cdots\exp(\mu_1 X_k)
\end{array}
\end{equation}
where $k=i$ or $j$, and $\ell=j$ or $i$, according to the parity of $m=3,4,6$;
$\mu_1,\ldots,\mu_m$ are linear functions in 
$a=\lambda_i,b=\lambda_j$ for $k=1,\ldots,m$. 
Denote the left hand side and the right hand side 
of (\ref{EE1}), regarded as functions 
of $\lambda=(\lambda_i)_{i\in I}$, by $\Phi(\lambda)$ and by $\Psi(\lambda)$, 
respectively. 
Since $\Phi(0)=\Psi(0)$ in each case,  equality (\ref{EE1}) reduces to proving
\begin{equation} 
\partial_{\lambda_i}(\Phi(\lambda))\,\Phi(\lambda)^{-1}=
\partial_{\lambda_i}(\Psi(\lambda))\,\Psi(\lambda)^{-1}\qquad(i\in I).
\end{equation}
For an arbitrary element $\Lambda\in L$, we define the derivation with respect to 
$\Lambda$ by
\begin{equation}
\partial_{\Lambda}=\sum_{i\in I}  \br{h_i,\Lambda} \partial_{\lambda_i}. 
\end{equation}
We will derive the formula 
\begin{equation}
\partial_{\Lambda}(\Phi(\lambda))\,\Phi(\lambda)^{-1}=
\partial_{\Lambda}(\Psi(\lambda))\,\Psi(\lambda)^{-1}\qquad(\Lambda\in L)
\end{equation}
from a general statement for operators defined as certain product
of canonical transformations of the form $R_i(t)=\exp(t\,X_i)$. 

\par\medskip
For the moment, we fix an arbitrary element $\Lambda\in L$. 
For each operator
$P=P(\lambda)\in\End_{\BC}(\widetilde{\CR})[[\lambda]]$ with 
invertible leading term $P(0)$, 
we set
\begin{equation}
L_P=\partial_\Lambda(P(\lambda)) P(\lambda)^{-1}. 
\end{equation}
Then, as to the composition of two operators, we have 
\begin{equation}
L_{PQ}=L_P+P(\lambda) L_Q P(\lambda)^{-1}. 
\end{equation}
If $P$ is expressed as 
\begin{equation}
P(\lambda)=\exp(a(\lambda) X),\quad X=\pad(f)
\end{equation}
with some $f\in \widetilde{\CR}$ and a linear function $a=a(\lambda)$, 
then we have 
\begin{equation}
\begin{array}{l}
\smallskip
L_P=\partial_\Lambda(a) X=\pad(\partial_\Lambda(a)f),\\
P(\lambda)L_QP(\lambda)^{-1}=\exp(a X)L_Q\exp(-aX)
=\exp(a\ad(X))(L_Q). 
\end{array}
\end{equation}
Suppose furthermore that $L_Q=\pad(g(\lambda))$ for some 
$g(\lambda)\in \widetilde{\CR}[[\lambda]]$.
Since $X=\pad(f)$, by the Jacobi identity for the Poisson bracket, we have
\begin{equation}
\exp(a\,\ad(X))(L_Q)=\exp(a\,\ad(X))(\pad(g(\lambda)))=\pad(\exp(a
X)g(\lambda)).
\end{equation}
Hence we obtain the following expression for $L_{PQ}$:
\begin{equation}\label{OB}
L_{PQ}=\pad(\partial_\Lambda(a) f+\exp(a X)g(\lambda)). 
\end{equation}
We make use of these observations to analyze the action of $\hW$.

\par\medskip
For each $w=s_{j_1}s_{j_2}\cdots s_{j_p}\in \hW$, we define the operator
$G_w(\lambda)$ by 
\begin{equation}\label{G} 
w=G_w(\lambda) r,\quad 
G_w(\lambda)=\exp(\mu_1 X_{j_1})\exp(\mu_2 X_{j_2})\cdots\exp(\mu_{p}X_{j_p}),
\end{equation}
where $r=\rho(w)=r_{j_1}\cdots r_{j_p}$ and $\mu_k=r_{j_1}\cdots
r_{j_{k-1}}(\lambda_{j_k})$ for 
$k=1,\ldots,p$.
In the following we use the multiplicative notation 
\begin{equation}
\log(\varphi_i ^a \varphi_j^b\cdots)
=a \log\varphi_i+b\log\varphi_j+\cdots\quad
(a,b,\cdots\in\BZ). 
\end{equation}
Also, for a power series in the form
\begin{equation}
f(\lambda)=f_0+\sum_{|\nu|>0}f_\nu \lambda^{\nu}\qquad(f_\nu\in \CR)
\end{equation}
(with the multi-index notation)
such that $f_0$ is a product of $\varphi_j$'s, we write 
\begin{equation}
\log f(\lambda)=\log f_0+ 
\log(1+f_0^{-1}\sum_{|\nu|>0}f_\nu \lambda^{\nu}),
\end{equation}
understanding the last term by its Taylor expansion. 

\WSCI{\par\medskip}
\begin{prop}\label{PropKey} For any element 
$w=s_{j_1}\cdots s_{j_p}\in \hW$ and 
$\Lambda\in L$, one has 
\begin{equation}
\partial_\Lambda(G_w(\lambda))G_w(\lambda)^{-1}
=-\pad(\log\phi_w(\rho(w)^{-1}.\Lambda)).
\end{equation}
\end{prop} 
\par\smallskip
\noindent{\sl Proof.}\quad
If $P=\exp(\lambda_{j}X_{j})$, 
\begin{equation}
L_P=\pad(\br{h_{j},\Lambda}\log\varphi_j)=\pad(\log \varphi_j^{\br{h_j,\Lambda}})
=-\pad(\log \varphi_j^{\br{h_j,r_j\Lambda}}).
\end{equation}
This proves the case where $p=1$. 
When $w=s_{j_1}\cdots s_{j_p}$ and $r=\rho(w)=r_{j_1}\cdots r_{j_p}$, 
by using the expression (\ref{G}), we set 
\begin{equation}
P=\exp(\lambda_{j_1} X_{j_1}),\quad Q=\exp(\mu_2 X_{j_2})\cdots \exp(\mu_p X_{j_p})
\end{equation}
so that $G_w=PQ$. 
By the induction hypothesis, we have
$$
\partial_\Lambda(G_{w'})
G_{w'}^{-1}=-\pad(\log\phi_{w'}({r'}^{-1}.\Lambda))
$$
for $w'=s_{j_1}w=s_{j_2}\cdots s_{j_p}$, $r'=r_{j_1}r=r_{j_2}\cdots r_{j_p}$,. 
Since $Q=r_{j_1}G_{w'} r_{j_1}$, we compute
$$
\partial_\Lambda(Q) Q^{-1}=r_{j_1}\partial_{r_{j_1}\Lambda}(G_{w'}){G_{w'}}^{-1}r_{j_1}
=-\pad(\log r_{j_1}(\phi_{w'}({r}^{-1}\Lambda))).
$$
Hence by (\ref{OB}) we obtain
\begin{equation}
\begin{array}{l}
\smallskip
\partial_\Lambda(G_w){G_w}^{-1}=L_{PQ}\\
\smallskip
=-\pad(\log \varphi_{j_1}^{\br{h_{j_1},r_{j_1}\Lambda}}+
\log \exp(\lambda_{j_1}X_{j_1})r_{j_1}(\phi_{w'}({r}^{-1}\Lambda)))\\
\smallskip
=-\pad(\log(\varphi_{j_1}^{\br{h_{j_1},r'r^{-1}\Lambda}}
s_{j_1}(\phi_{w'}(r^{-1}\Lambda))))\\ 
\smallskip
=-\pad(\log \phi_w(r^{-1}.\Lambda))
\end{array}
\end{equation}
for any $\Lambda\in L$. \  \QED
\par\medskip

When $(a_{ij},a_{ji})=(-k,-1)$ for  $k=1,2,3$, we set 
\begin{equation}
\Phi(\lambda)=G_{s_is_js_i\cdots}(\lambda), \quad \Psi(\lambda)=G_{s_js_is_j\cdots}(\lambda). 
\end{equation}
Then Proposition \ref{PropKey} implies
\begin{equation}
\begin{array}{l}
\smallskip
\partial_\Lambda(\Phi(\lambda))\Phi(\lambda)^{-1}
=-\pad(\log \phi_{s_is_js_i\cdots}(r^{-1}.\Lambda))\\
\partial_\Lambda(\Psi(\lambda))\Psi(\lambda)^{-1}
=-\pad(\log \phi_{s_js_is_j\cdots}(r^{-1}.\Lambda)),
\end{array}
\end{equation}
where $r=r_ir_jr_i\cdots=r_jr_ir_j\cdots$. 
Since $\phi_{s_is_js_i\cdots}=\phi_{s_js_is_j\cdots}$ by Lemma \ref{LemB}, we see
\begin{equation}
\partial_\Lambda(\Phi(\lambda))\Phi(\lambda)^{-1}
=\partial_\Lambda(\Psi(\lambda))\Psi(\lambda)^{-1}
\end{equation}
for any $\Lambda\in L$.  
Since $\Phi(0)=\Psi(0)$, we have $\Phi(\lambda)=\Psi(\lambda)$, namely
\begin{equation}
G_{s_is_js_i\cdots}(\lambda)=G_{s_js_is_j\cdots}(\lambda),
\end{equation}
which proves the Yang-Baxter equation of (\ref{YBE}).  
Since $s_is_js_i\cdots=G_{s_is_js_i\cdots}r$ and 
$s_js_is_j\cdots=G_{s_js_is_j\cdots}r$, we obtain the braid relation 
of (\ref{BR}) on $\CR[[\lambda]]$ as desired.  
Since we have already seen the validity of braid relation on 
$\tau$-functions in Corollary \ref{cor-BRtau}, 
this completes the proof of Theorems \ref{thm-E}. 

\rem{
In the formulation of Section 1, we assumed for simplicity 
that $\CA_0$ is {\em generated as a Poisson algebra}, 
by a set of nonzero elements $\varphi_i$ ($i\in I$) satisfying the 
Serre relations.  
As can be seen from the proof of this section, 
for the validity of Theorems \ref{thm-A} and \ref{thm-B}, 
we have only to 
assume that the following three conditions are satisfied: 
\par\medskip
\begin{tabular}{cl}
(i) & $\CA_0$ has no zerodivisors.\cr
(ii) & $\CA_0$ contains a set of nonzero elements $\varphi_i$ ($i\in I$) 
such that \cr 
& $\pad(\varphi_i)^{-a_{ij}+1}\varphi_j=0$ ($i\ne j$). \cr
(iii) & The action of $\pad(\varphi_i)$ on $\CA_0$ is {\em locally nilpotent}
for each $i\in I$. 
\end{tabular}
}

\section{Lie theoretic background}

In this section, we assume that the GCM $A=(a_{ij})_{ij\in I}$ is 
{\em symmetrizable} and that the indexing set $I$ is finite. 
We fix a set of nonzero rational numbers $\epsilon_{i}$ ($i\in I$) such 
that $a_{ij}\epsilon_j=a_{ji}\epsilon_i $ for $i,j\in I$. 
We will explain below a geometric background 
of our realization of the Weyl group $W=W(A)$ 
in the language of Kac-Moody Lie algebras and Kac-Moody groups. 


\subsection{Birational action of $\,W$ on the Borel subalgebra}

As to Kac-Moody Lie algebras, 
we will basically follow the notation of Kac\WSCI{.\cite{K} }\PREP{ \cite{K}. }
Let $\frg=\frg(A)$ be the Kac-Moody Lie algebra associated with 
the symmetrizable generalized Cartan matrix $A$, and consider the 
triangular decomposition $\frg=\frn\oplus \frh \oplus \frn_{-}$,
where $\frn$ (resp. $\frn_-$) is the nilpotent Lie subalgebras of 
$\frg$ generated by $e_i$ (resp. $f_i$).  
We denote by $\frg=\bigoplus_{n\in\BZ} \frg_n$ the gradation 
of $\frg$ with respect to the degree (height) such that 
$\deg e_i=1$, 
$\deg f_i=-1$ ($i\in I$) and $\deg h=0$ ($h\in \frh$).
Note that each $\frg_{n}$ is finite dimensional and that 
$\frn=\oplus_{n>0}\, \frg_n$, 
$\frh=\frg_0$, $\frn_{-}=\oplus_{n<0}\, \frg_{n}$. 
We denote by 
$\frg=\frh\oplus\sum_{\alpha\in\Delta}\frg_{\alpha}$ 
the root space decomposition of $\frg$ with respect to the 
Cartan subalgebra $\frh$, so that 
$\frg_n=\sum_{\mbox{\small ht}(\alpha)=n} \frg_\alpha$ 
($n\in \BZ$, $n\ne 0$). 

\par\medskip
The symmetrizable Kac-Moody Lie algebra $\frg=\frg(A)$ has an 
$\ad$-invariant symmetric bilinear form 
$(\ |\ ) : \frg\times\frg\to\BC$ which induces nondegenerate 
pairings $\frh\times\frh\to\BC$ and 
$\frg_{-\alpha}\times\frg_{\alpha}\to \BC$  ($\alpha\in\Delta$).
We take the normalization of $(\ |\ )$ such that
\begin{equation}
(h_i|h_j)=a_{ij}\epsilon_j, \quad 
(e_i|f_j)=\delta_{ij}\epsilon_j\quad(i,j\in I).
\end{equation} 
Through this symmetric bilinear form, we have an isomorphism 
$\iota :\frg\iso \frg^\ast $ of $\BC$-vector spaces such 
that $\iota(y)(x)=(y|x)$ for $y, x\in \frg$,
where $\frg^\ast=\oplus_{n\in\BZ} \frg_n^\ast$ is the 
{\em graded dual} of $\frg$. 
We define the bracket 
$\pbr{\,,\ } : \frg^\ast\times\frg^\ast\to\frg^\ast$
by
\begin{equation}
\pbr{\iota(y),\iota(y')}=-\iota([y,y'])\qquad(y,y'\in\frg). 
\end{equation}
This bracket naturally extends to a Poisson bracket of 
$\Sym(\frg^\ast)$.  
Note that
each element of the Poisson algebra $\Sym(\frg^\ast)$ 
is regarded naturally as a polynomial function on $\frg$. 
If $\varphi=\iota(y)$ and $\psi=\iota(y')$ ($y,y'\in\frg$),
the $\ad$-invariance of $(\ |\ )$ implies 
\begin{equation}\label{eq-ad}
(\pad(\varphi)\psi)(x)
=-(\ad(y)y'|x)=(y'|\ad(y)x)=\psi(\ad(y)x) 
\end{equation}
for any $x\in\frg$. 
Hence we have

\WSCI{\par\medskip}
\begin{lem}\label{lem-adad}
Let $y\in\frg$ and set $\varphi=\iota(y)\in \frg^\ast$.  
Then, for any $\psi\in\Sym(\frg^\ast)$, one has
\begin{equation}\label{eq-exp}
(\exp(t\,\pad(\varphi))\psi)(x)=\psi(\exp(t\,\ad(y)) x)\qquad(x\in\frg),
\end{equation}
where $t$ is a formal parameter. 
\end{lem}
\WSCI{\par\medskip}
\noindent
One has only to verify (\ref{eq-exp}) for $\psi\in \frg^\ast$, since 
the both sides of (\ref{eq-exp}) is multiplicative with respect to $\psi$. 
This equality for $\psi=\iota(y')$ ($y'\in\frg$) follows 
from (\ref{eq-ad}).
\par\medskip
We consider now the formal completion 
$\widehat{\frn}=\projlim{k} \frn/\ad(\frn)^k\frn$
of the nilpotent subalgebra $\frn\subset\frg$
with respect to the height, and set 
$\widehat{\frb}=\frh\oplus\widehat{\frn}$. 
The coordinate rings of $\frh$, $\widehat{\frn}$ 
and $\widehat{\frb}$ are defined by 
\begin{equation}
\CO(\frh)=\Sym(\frh^\ast), \quad 
\CO(\widehat{\frn})=\Sym(\frn^\ast)\quad\mbox{and}\quad
\CO(\widehat{\frb})=\Sym(\frh^\ast\oplus\frn^\ast), 
\end{equation}
respectively, 
where $\frn^\ast=\oplus_{n=1}^\infty \frg_n^\ast$ stands for 
the graded dual of $\frn$. 
Note that $\iota: \frg\iso\frg^\ast$ induces the 
isomorphism $\frn_{-}\iso\frn^{\ast}$
and that $\CO(\widehat{\frn})=\Sym(\frn^\ast)$ can be 
regarded as a Poisson subalgebra of $\Sym(\frg^\ast)$. 
Let $\pbr{\varphi_i}_{i\in I}$ be the basis of 
$\frg_1^\ast$ dual to the basis $\pbr{e_i}_{i\in I}$ of 
$\frg_1$: $\varphi_i(e_j)=\delta_{ij}$ ($i,j\in I$). 
Then 
$\CO(\widehat{\frn})$ is generated as a Poisson algebra 
by $\pbr{\varphi_i}_{i\in I}$, and these generators 
satisfy the Serre relations
\begin{equation}
\pad(\varphi_i)^{-a_{ij}+1}(\varphi_j)=0\qquad(i\ne j)
\end{equation}
with respect to the Poisson bracket, since 
$\varphi_i=\iota(f_i)/\epsilon_i$ ($i\in I$).
Note also that the simple roots $\alpha_i$ ($i\in I$) are regarded 
as functions on $\widehat{\frb}$ of height 0.
Hence the coordinate ring $\CO(\widehat{\frb})$ contains 
the subalgebra $\BC[\alpha]\otimes\CO(\widehat{\frn})$. 

In this setting we explain a procedure to define a birational action 
of the Weyl group $W(A)$ on the Borel subalgebra $\widehat{\frb}$. 
For each $i$, consider the adjoint action of the 
1-parameter subgroup $\exp(tf_i)$ ($t\in\BC$) on $\widehat{\frb}$. 
Note that $\exp(t\,\ad(f_i))$ is well defined as a mapping 
$\widehat{\frb}\to \BC f_i\oplus \widehat{\frb}$, 
since, for each $\beta\in\Delta_{+}$, 
there are only a finite number of positive roots 
of the form $\beta+n\alpha_i$ ($n=0,1,2,\ldots$).
For an element $x\in\widehat{\frb}$ given, 
we try to find a value of $t$ 
such that 
\begin{equation} 
\Ad(\exp(t f_i))x=\exp(t\,\ad(f_i))x \in \widehat{\frb}. 
\end{equation}
A general element $x\in \widehat{\frb}$ can be expressed 
in the form
\begin{equation}
x=h+\sum_{j\in I} c_j e_j+(\mbox{components of height $\ge 2$}),
\end{equation}
where $h\in\frh$ and  $c_j=\varphi_j(x)$  ($j\in I$).
By a simple computation, we see that 
$\exp(t\,\ad(f_i))x$ takes the form 
\begin{equation}
t(\br{h,\alpha_i}-t\,c_i)f_i+(h-t\,c_i h_i)+(\mbox{components of height $\ge1$}).
\end{equation}
Hence we have a nontrivial solution
\begin{equation}
t=\frac{\br{h,\alpha_i}}{c_i}=\frac{\alpha_i(x)}{\varphi_i(x)},
\end{equation}
provided that $c_i=\varphi_i(x)\ne 0$. 
Note that, when $t=\br{h,\alpha_i}/c_i$, one has
\begin{equation}\label{eq-adx}
\exp(t\,\ad(f_i))x=h-\br{h,\alpha_i}h_i+(\mbox{components of height $\ge1$}).
\end{equation}
From this argument, for each $i\in I$, we obtain a regular mapping
\begin{equation}
\sigma_i :\ x\mapsto \Ad(\exp(t_i(x)f_i))x\ :\  
\widehat{\frb}\setminus\pbr{\varphi_i=0}\to \widehat{\frb},
\end{equation}
where $t_i=\alpha_i/\varphi_i$. 
Passing to the coordinate ring, we denote by 
$\sigma_i^\ast$ the automorphism of 
the field $\CK(\widehat{\frb})$ of rational functions on $\widehat{\frb}$,
corresponding to the birational mapping 
$\sigma_i: \widehat{\frb}\cdots\!\to\widehat{\frb}$.

\WSCI{\par\medskip}
\begin{prop}\label{prop-sigma}
For each $i\in I$, 
the automorphism $\sigma_i^\ast$ of $\CK(\widehat{\frb})$ is 
characterized by the following two properties: 
\begin{equation}
\begin{array}{ll}
\smallskip
(1)\qquad\sigma_i^\ast(\lambda)=\lambda-\alpha_i\br{h_i,\lambda}
&(\lambda\in\frh^\ast),\cr
(2)\qquad
\sigma_i^\ast(\psi)
=\exp\left(\dfrac{\lambda_i}{\varphi_i}\pad(\varphi_i)\right)\,\psi\quad
&(\psi\in\CO(\widehat{\frn})),
\end{array}
\end{equation}
where $\lambda_i=\epsilon_i\alpha_i$. 
\end{prop}
\WSCI{\par\medskip}
\noindent
The first equality for $\lambda\in\frh^\ast$ is a consequence of 
the expression (\ref{eq-adx}).
If $\psi\in \CO(\widehat{\frn})$, by Lemma \ref{lem-adad} one has
\begin{equation}
\begin{array}{l}
\smallskip
\sigma_i^\ast(\psi)(x)=\psi(\exp(t_i(x)\ad(f_i))x)\cr
\smallskip
\phantom{\sigma_i^\ast(\psi)(x)}
=(\exp(\epsilon_i \,t_i(x)\pad(\varphi_i))\psi)(x)\cr
\phantom{\sigma_i^\ast(\psi)(x)}
=(\exp(\epsilon_i \,t_i\,\pad(\varphi_i))\psi)(x)
\end{array}
\end{equation}
for any $x\in \widehat{\frb}\setminus\pbr{\varphi_i=0}$.
This proves the second statement, 
since $\epsilon_i\,t_i=\epsilon_i\,\alpha_i/\varphi_i$. 

\par\medskip
These automorphisms $\sigma_i^\ast$ ($i\in I$) are 
essentially the same as $s_i$ ($i\in I$) in Theorem \ref{thm-A}
for the nilpotent Poisson algebra $\CA_0=\CO(\widehat{\frn})$. 
In fact, 
the polynomial ring $\BC[\lambda]$ with 
$\lambda_i=\epsilon_i\alpha_i$ is a subalgebra of 
$\CO(\frh^\ast)$, and 
each automorphism $\sigma_i^\ast$ coincides
with $s_i$ on the subfield 
$\CK=Q(\BC[\lambda]\otimes\CO(\widehat{\frn}))$ of 
$\CK(\widehat{\frb})$.

Although the construction of $\sigma_i^\ast$ gives an origin of our 
birational automorphisms $s_i$, 
it does not explain why these automorphisms satisfy 
the fundamental relations for the generators of the 
Weyl group $W=W(A)$. 
In the following subsections, we will show that the birational 
automorphism $\sigma_i$ of $\widehat{\frb}$ arise 
naturally from the birational {\em dressing\/} action 
of a lift of the Weyl group on the Borel subgroup, induced 
through the Gauss decomposition
in the Kac-Moody group. 
This construction also provides an alternative 
proof of Theorem \ref{thm-A} 
for the case when the GCM $A$ is symmetrizable and 
$\CA_0=\CO(\widehat{\frn})$.

\subsection{Birational action of $\dot{W}$ on the Borel subgroup}

We denote by $\frg'=\frn\oplus \frh'\oplus \frn_{-}$ be the 
subalgebra of $\frg$, generated by the Chevalley generators 
$e_i, h_i, f_i$ 
($i\in I$), where $\frh'=\bigoplus_{i\in I} \BC \,h_i$. 
In the following we will make use of a variant of Kac-Moody 
group $G$ for the Lie algebra $\frg'$.
Although there are several variants of Kac-Moody 
groups\WSCI{,\cite{KP1,K,M,S}}
\PREP{(\cite{KP1}, \cite{K}, \cite{M}, \cite{S}),} 
we will work now with the one as in Slodowy\WSCI{.\cite{S} }\PREP{\cite{S}.} 
(In the argument below,  we can also use the infinite 
dimensional scheme of Kashiwara\WSCI{,\cite{K} }\PREP{ \cite{K},}
since we do not 
need the group structure on the whole $G$.)

We first recall some properties of $G$ which will be used 
in our argument. 
\par\smallskip\noindent
(1) \ The group $G$ contains the Borel subgroup
\begin{equation}
B=T\ltimes U,\quad T=\Spec(\BC[L]),\quad U=\exp(\widehat{\frn})
\end{equation}
where $L=\Hom_{\,\BZ}(Q^\vee,\BZ)$, and 
$\widehat{\frn}=\projlim{k}
\frn/\ad(\frn)^k(\frn)$ is the completion of $\frn$ with respect 
to the height.
\par\noindent
(2) \ For each $i\in I$, $G$ contains a parabolic subgroup $P_i\supset B$
associated with the simple root $\alpha_i$.  
It has the Levi 
decomposition 
\begin{equation}
P_i=G_i\ltimes U_i, \quad U_i=\exp(\widehat{\frn}_i)
\end{equation}
where $G_i$ the group of the Lie algebra 
$\frg_i'=\BC e_i\oplus \frh'\oplus \BC f_i$ with the Cartan subgroup $T$, 
and 
$\widehat{\frn}_i=\projlim{k}\frn_i/\ad(\frn)^k(\frn_i)$ 
is the completion of 
$\frn_i=\bigoplus_{\alpha\in\Delta_+\setminus\pbr{\alpha_i}}\
\frg_{\alpha}$ 
with respect to the height. 
\par\noindent
(3) \ The group $G$ has a subgroup $U_{-}$ such that $U_{-}\cap B=\pbr{1}$ 
and that all the 1-parameter subgroups associated with the negative real roots
are contained in $U_{-}$. 
\par\noindent
(4)  \ For each dominant integral weight $\Lambda\in P_{+}$, 
(i.e., $\Lambda\in\frh^\ast$ such that 
$\br{h_i,\Lambda}\in \BZ_{\ge0}$ for all $i\in I$), 
the irreducible $\frg$-module $L(\Lambda)$ 
with highest weight $\Lambda$ can be integrated 
to a $G$-module. 
\par\medskip
For each $i\in I$, we set
\begin{equation}
\ds_i=\exp(-e_i)\exp(f_i)\exp(-e_i)\in G_i\subset P_i,
\end{equation}
and denote by $\dot{W}$ the subgroup of $G$ generated 
by $\ds_i$ ($i \in I$).  
These elements satisfy $\ds_j^4=1$ and 
the braid relations
\begin{equation}
\ds_i\ds_j\ds_i\cdots=\ds_j\ds_i\ds_j\cdots \quad
(\mbox{$m_{ij}$ factors on each side})
\end{equation}
in $G$, for any $i,j\in I$ with $a_{ij}a_{ji}=0,1,2,3$. 
In what follows, 
we denote by $W=W(A)$ the Weyl group on the generators $s_i$ 
$(i\in I)$.  
We have a natural group homomorphism $\dot{W}\to W$ which 
maps $\dot{s}_i$ to $s_i$ ($i\in I$). 

\par\medskip
We now express a general element of the Borel subgroup $B=T\ltimes U$ as 
\begin{equation}
b=\tau \exp(\xi), \quad \tau=\prod_{i\in I} {\tau_i}^{h_i}\in T,\quad
\xi=\sum_{n=1}^\infty\,\xi_n\in \widehat{\frn}\ \ \ (\xi_n\in \frg_n).
\end{equation}
In this decomposition, the coordinate ring $\CO(U)$ of $U$ is defined 
as the symmetric algebra $\Sym(\frn^{\ast})$ 
of the graded dual
of $\frn$.  
For a general element $\xi\in \widehat{\frn}$, 
we express the component $\xi_1$ of height 1 as 
$\xi_1=\sum_{i\in I}\psi_i e_i$, 
regarding $\psi_i$ ($i\in I$) as the part of coordinates of $U$ in height 1. 
These $\psi_i$ ($i\in I$) form a $\BC$-basis of $\frg_1^\ast$ dual 
to the $\BC$-basis $e_i$ ($i\in I$) of $\frg_1$.
Similarly we regard $\tau_i$ ($i \in I$) as the coordinates of the 
Cartan subgroup $T\simeq (\BC^{\times})^I$, 
and use the notation
$\tau^\Lambda=\prod_{i\in I}\,\tau_i^{\br{h_i,\Lambda}}
\in \CO(T)$
for any $\Lambda\in L$. 
Note that the coordinate ring $\CO(B)$ of $B$ is represented as 
the ring of Laurent polynomials $\CO(U)[\tau^{\pm1}]$ in $\tau_i$ ($i\in I$) 
with coefficients in $\CO(U)=\Sym(\frn^\ast)$. 

Let $\dot{w}$ be an arbitrary element of $\dot{W}$. 
Suppose that $\exp(\xi)\in U$ is {\em generic} 
so that $\exp(\xi)\,\dot{w}\in U_{-}\,B$. 
Then, for  $b=\tau\exp(\xi)\in B$ with an arbitrary $\tau\in T$, 
the product $b\,\dot{w}$ can be decomposed uniquely in the form
\begin{equation}
b\,\dot{w}= g^{-1}\,\widetilde{b},\quad
g\in U_{-},\quad 
\widetilde{b}\in B. 
\end{equation}
This type of {\em Gauss decomposition} 
(Birkhoff decomposition) is determined 
by carrying out the Gauss decomposition for the right action 
of $\dot{s_i}^{\pm1}$ ($i\in I$) repeatedly. 
For the argument below, we examine the case of $\dot{w}=\dot{s}_i$ 
in some detail. 

\WSCI{\par\medskip}
\begin{prop}\label{prop-GD1}
For each $i\in I$, any element $b=\tau\exp(\xi)$ with $\psi_i\ne 0$ 
is decomposed uniquely in the form $b\,\dot{s}_i=g_i^{-1}\,\widetilde{b}$ 
with $g_i\in U_{-}$ and $\widetilde{b}=\widetilde{\tau}\exp(\widetilde{\xi})\in B$.  
The element $g$ and $\widetilde{\tau}$ are given by
\begin{equation}
g_i=\exp\left(\frac{-1}{\psi_i\tau^{\alpha_i}}f_i\right), 
\quad
\widetilde{\tau}=\tau\,\psi_i^{h_i},
\end{equation}
where $\tau^{\alpha_i}=\prod_{k\in I} \tau_k^{\br{h_k,\alpha_i}}
=\prod_{k\in I} \tau_k^{a_{ki}}$. 
Furthermore, the correspondence $\xi\to\widetilde{\xi}$ is 
represented by an algebra homomorphism 
$\CO(U)\to \CO(U)[\psi_{i}^{-1}]$. 
\end{prop}
\WSCI{\par\medskip}
\noindent
The Gauss decomposition of $b\,\dot{s}_i$ is determined by employing 
the Levi decomposition $P_i=G_i\ltimes U_{i}$ and the following lemma
in $G_i$. 

\WSCI{\par\medskip}
\begin{lem} \label{lem-GD1}
Let $a,b\in\BC$ and $c\in \BC^{\times}$. 
For each $i\in I$, one has the identity
\begin{equation}
\exp(a\,e_i)\, c^{h_i}\exp(b\,f_i)=
\exp(y\,f_i)\, z^{h_i}\exp(x\,e_i)
\end{equation}
in $G_j$ unless $ab+c^2=0$,  where 
\begin{equation}
x=\frac{a}{ab+c^2},\quad z=\frac{ab+c^2}{c},\quad y=\frac{b}{ab+c^2}.
\end{equation}
\end{lem}
\WSCI{\par\medskip}
\noindent
{\em Proof of Proposition \ref{prop-GD1}.}\quad
By decomposing $b=\tau\exp(\xi)\in B$ in the form 
\begin{equation}
b=\tau\,\exp(\psi_i e_i)\exp(\xi'),\quad \xi'\in \widehat{\frn}_i, 
\end{equation}
we get
\begin{equation}
b\,\dot{s}_i=\tau\,\exp(\psi_i e_i)\exp(\xi') \dot{s}_i
=\tau \exp(\psi_i e_i) \dot{s}_i
\exp(\xi''),
\end{equation}
where 
$\xi''=\Ad(\dot{s}_i^{-1})(\xi')$.
Note that $\exp(-\ad(f_i)) : \widehat{\frn}_i\to\widehat{\frn}_i$ is 
well defined since, for any $\beta\in \Delta_{+}\setminus\pbr{\alpha_i}$, 
there are only a finite number of positive roots in the form 
$\beta+n\alpha_i$ ($n\ge 0$). 
Then by Lemma \ref{lem-GD1} we have
\begin{equation}
\begin{array}{l}
\smallskip
\exp(\psi_i e_i)\dot{s}_i=\exp(\psi_i e_i)
\exp(-e_i)\exp(f_i)\exp(-e_i)\cr
\smallskip
\phantom{\exp(\psi_i e_i)\dot{s}_i}
=\exp((\psi_i-1)e_i)\exp(f_i)\exp(-e_i)\cr
\smallskip
\phantom{\exp(\psi_i e_i)\dot{s}_i}
=\exp(\frac{1}{\psi_i}f_i)\psi_i^{h_i}
\exp((1-\frac{1}{\psi_i})e_i)\exp(-e_i)\cr
\phantom{\exp(\psi_i e_i)\dot{s}_i}
=\exp(\frac{1}{\psi_i}f_i)\psi_i^{h_i}\exp(\frac{-1}{\psi_i}e_i). 
\end{array}
\end{equation}
Hence,
\begin{equation}
\begin{array}{l}
\smallskip
b\,\dot{s}_i=\tau\exp(\frac{1}{\psi_i}f_i)
\psi_i^{h_i}\exp(\frac{-1}{\psi_i}e_i)\exp(\xi'')\cr
\phantom{b\,\dot{s}_i}
=\exp(\frac{1}{\psi_i\tau^{\alpha_i}}f_i) \ \tau \,\psi_i^{h_i}\,
\exp(\widetilde{\xi}),
\end{array}
\end{equation}
where $\exp(\widetilde{\xi})=\exp(\frac{-1}{\psi_i}e_i)\exp(\xi'')$.
This computation implies that $\widetilde{\xi}$ is regular on 
$\psi_i\ne 0$, and $g=\exp(\frac{-1}{\psi_i\tau^{\alpha_i}}f_i)$,
$\widetilde{\tau}=\tau\, \psi_i^{h_i}$, as desired. 
\QED
\par\medskip 
We now consider the Gauss decomposition of  
$b \,\dot{w}$ for an arbitrary $\dot{w}\in\dot{W}$: 
\begin{equation}
b\,\dot{w} =g^{-1}\,\widetilde{b},\quad g\in U_{-},\quad \widetilde{b}\in B,
\end{equation}
for generic $b\in B$. 
We denote the two elements $g$ and 
$\widetilde{b}$ above 
by $g(\dot{w},b)\in U_{-}$ and $b\ast\dot{w}\in B$,
respectively. 
From this Gauss decomposition of $b\,\dot{w}$, we obtain a 
mapping $\ast\dot{w}$ 
from an open dense subset of $B$, 
which sends generic $b\in B$ to 
$b\ast\dot{w}=g(\dot{w},b)\,b\,\dot{w}\in B$.  
Note that the uniqueness of Gauss decomposition implies the 
following $1$-cocycle relation for $g(\dot{w},b)$:
\begin{equation}\label{eq-gcocycle}
g(1,b)=1,\quad 
g(\dot{w}_1\dot{w}_2,b)=
g(\dot{w}_2,b\ast\dot{w}_1)\, g(w_1,b)
\end{equation}
for any $\dot{w}\in \dot{W}$ and for generic $b\in B$. 
From Proposition \ref{prop-GD1}, we see that, for each $i\in I$, 
$g_i=g(\dot{s}_i,b)$ is defined for $b=\tau\exp(\xi)$ with $\psi_i\ne 0$ by 
\begin{equation}
g_i=g(\dot{s}_i,b)=\exp\left(\frac{-1}{\psi_i\,\tau^{\alpha_j}} f_i\right ).
\end{equation}
The element $g(\dot{w},b)$ for general $\dot{w}\in\dot{W}$ 
is determined inductively from 
$g(\dot{s}_j,b)$ by the cocycle relation (\ref{eq-gcocycle}). 
Also, we obtain a birational right action of $\dot{W}$ on $B$:
\begin{equation}
b\ast 1=b,\quad 
b\ast(\dot{w}_1\dot{w_2})=(b\ast\dot{w}_1)\ast\dot{w}_2
\end{equation}
for any $\dot{w}_1,\dot{w}_2\in \dot{W}$ and 
for generic $b\in B$. 
In terms of the coordinate rings, this birational mapping 
$\ast \dot{w}: B\cdots\!\to B$ 
defines an algebra  automorphism
\begin{equation}
R_{\dot{w}}: \CK(U)[\tau^{\pm}] \to \CK(U)[\tau^{\pm}],
\end{equation}
$\CK(U)$ being the field of rational functions on $U$, 
such that
\begin{equation}
(R_{\dot{w}}\psi)(b)=\psi(b\ast\dot{w})=\psi(g(\dot{w},b)\, b\, \dot{w})
\end{equation}
for any $\psi\in\CK(U)[\tau^{\pm1}]$ and for generic $b\in B$. 
Note also that 
\begin{equation}
R_{1}=1,\quad R_{\dot{w}_1\dot{w}_2}=R_{\dot{w}_1}R_{\dot{w}_2}
\quad(\dot{w}_1,\dot{w}_2\in \dot{W}). 
\end{equation} 
These $R_{\dot{w}}$ give a realization of $\dot{W}$ as a group of 
automorphisms of the ring $\CK(U)[\tau^{\pm1}]$ of Laurent polynomials. 
We remark that, if $\dot{w}=\dot{s}_{j_1}\cdots\dot{s}_{j_p}$, 
the homomorphism 
$R_{\dot{w}}$ is determined as the composition 
$R_{\dot{w}}=R_{\dot{s}_{j_1}}\cdots R_{\dot{s}_{j_p}}$.

\subsection{Passage to the Borel subalgebra} 
We now consider to transfer the birational action of $\dot{W}$ on 
the Borel subgroup $B$ to that on the Borel subalgebra $\widehat{\frb}$
by means of the adjoint action. 

We say that an element $h$ of the Cartan subalgebra $\frh$ is 
{\em regular} if $\br{h,\alpha}\ne 0$ for any positive root
$\alpha\in \Delta_{+}$.  
We denote 
the set of all regular elements in $\frh$ by 
$\frh_{\mbox{\scriptsize reg}}$, and set 
$\widehat{\frb}_{\mbox{\scriptsize rs}}
=\frh_{\mbox{\scriptsize reg}}+\widehat{\frn}$, 
so that 
$\CO(\widehat{\frb}_{\mbox{\scriptsize rs}})=
\CO(\frh_{\mbox{\scriptsize reg}})\otimes \CO(\widehat{\frn})$, 
$\CO(\frh_{\mbox{\scriptsize reg}})=
\Sym(\frh^\ast)[\alpha^{-1} \,(\alpha\in\Delta_{+})]$.

\WSCI{\par\medskip}
\begin{prop} \label{prop-Ad}
The mapping 
\begin{equation}
\mu: B\times \frh_{\mbox{\scriptsize\rm reg}} \to 
T\times \widehat{\frb}_{\mbox{\scriptsize\rm rs}}\ :\  
(b=\tau\exp(\xi),h)\mapsto (\tau,\Ad(b)(h))
\end{equation}
gives an isomorphism of affine spaces. 
\end{prop}
\par\medskip\noindent
{\em Proof.}\quad
Under the condition 
that $h\in\frh$ is regular, 
for $\tau\in T$ and 
$x\in \widehat{\frn}$ given, 
one can find a unique $\xi\in\widehat{\frn}$  
such that $\Ad(\tau\exp(\xi))h=h+x$, i.e., 
\begin{equation}
\exp(\ad(\xi))(h)=h+\Ad(\tau^{-1})(x)
\end{equation}
inductively with respect to the height. 
In fact, 
by the root space decomposition
\begin{equation}
\xi=\sum_{\alpha\in \Delta_{+}} \xi_\alpha,
\quad
x=\sum_{\alpha\in \Delta_{+}} x_\alpha\qquad
(\xi_{\alpha},x_\alpha\in \frg_{\alpha}), 
\end{equation}
the equation above is equivalently rewritten into the recurrence formulas
\begin{equation}\label{eq-rec}
\br{h,\alpha} \xi_{\alpha}+
\sum_{k=2}^{\mbox{\small ht}(\alpha)}\frac{1}{k!}
\sum_{\beta_1+\cdots+\beta_k=\alpha}
\br{h,\beta_k} 
\ad(\xi_{\beta_1})\cdots\ad(\xi_{\beta_{k-1}})(\xi_{\beta_{k}})
=-\tau^{-\alpha}x_\alpha
\end{equation}
for $\alpha\in \Delta_{+}$. 
In terms of the coordinate rings, these recurrence formulas
imply that the mapping $(b,h)\to \Ad(b)(h)$ is represented 
the $\CO(\frh)$-algebra homomorphism 
$\CO(\widehat{\frb})=\CO(\widehat{\frn})\otimes\CO(\frh)
\to \CO(B)\otimes\CO(\frh)$ and that 
it induces an isomorphism 
\begin{equation}
\begin{array}{l}
\smallskip
\mu^\ast : \ \ \CO(T)\otimes
\CO(\widehat{\frb}_{\mbox{\scriptsize\rm rs}})=
\CO(T)\otimes
\CO(\frh_{\mbox{\scriptsize\rm reg}})\otimes\CO(\widehat{\frn})\cr
\quad\qquad\iso \CO(B)\otimes\CO(\frh_{\mbox{\scriptsize\rm reg}})=
\CO(T)\otimes\CO(\frh_{\mbox{\scriptsize\rm reg}})\otimes\CO(U)\cr
\end{array}
\end{equation}
of $\CO(T)\otimes\CO(\frh_{\mbox{\scriptsize\rm reg}})$-algebras. 
This proves Proposition. \QED
\par\medskip
\noindent
The recurrence formulas (\ref{eq-rec}) imply the inclusions
\begin{equation}
\begin{array}{c}
\quad\mu^\ast (\CO(\widehat{\frn}))\subset\CO(\frh)\otimes 
\CO(U)\otimes \CO(T), \cr
(\mu^\ast)^{-1}(\CO(U))\subset \CO(\frh_{\mbox{\scriptsize reg}})\otimes 
\CO(\widehat{\frn})\otimes \CO(T).
\end{array}
\end{equation}
We remark that, as a special case of (\ref{eq-rec}) when 
$\alpha=\alpha_i$,  we have
\begin{equation}
\br{h,\alpha_i} \psi_i =-\tau^{-\alpha_i} \varphi_i \quad (h\in \frh,\ i\in I). 
\end{equation}
It means that, under the identification of Proposition \ref{prop-Ad}, 
we have 
\begin{equation}\label{eq-phipsi}
\psi_i=-\dfrac{\varphi_i}{\alpha_i}\tau^{-\alpha_i},\quad
\varphi_i=-\alpha_i \psi_i\,\tau^{\alpha_i},\quad
\mbox{where}\ \ \tau^{\alpha_i}=\prod_{k\in I}\tau_k^{a_{ki}}.
\end{equation}
(In the expression $\tau^{\alpha_i}$, 
$\alpha_i$ is treated as a simple root.  
Otherwise, $\alpha_i$ is regarded as a function on $\frh$.) 
\par\medskip
We now consider to transfer the birational action of $\dot{W}$ 
from $B$ to $\widehat{\frb}_{\mbox{\scriptsize rs}}$ by 
the isomorphism of Proposition \ref{prop-Ad}.  
The product space $B\times \frh_{\mbox{\scriptsize reg}}$ admits 
the right birational action of $\dot{W}$ defined by
\begin{equation}
(b, h)\ast\dot{w}=(b\ast \dot{w}, \,\Ad(\dot{w})^{-1}h)
=(g(\dot{w},b)\,b\,\dot{w},\Ad(\dot{w})^{-1}h)
\end{equation}
for generic $b\in B$ and for $h\in \frh_{\mbox{\scriptsize reg}}$.
We denote by the same symbol $\ast \dot{w}$ 
the birational action of $\dot{W}$ on 
$T\times \widehat{b}_{\mbox{\scriptsize rs}}$ obtained from that 
$B\times \frh_{\mbox{\scriptsize reg}}$ through the 
isomorphism 
$\mu : B\times \frh_{\mbox{\scriptsize reg}}\iso 
T\times \widehat{\frb}_{\mbox{\scriptsize rs}}$
of Proposition \ref{prop-Ad}. 
Note that $\mu$ is equivariant 
with respect to the right action of the torus $T$ 
on the first component of 
$B\times \frh_{\mbox{\scriptsize reg}}$ and  
$T\times \widehat{\frb}_{\mbox{\scriptsize rs}}$,
respectively.
Since $\dot{w}$ normalizes $T$, the birational 
action of $\dot{w}$ on 
$T\times \widehat{\frb}_{\mbox{\scriptsize rs}}$
passes to $\widehat{\frb}_{\mbox{\scriptsize rs}}$
through the second projection.
(For 
$(\tau,x)\in T\times \widehat{\frb}_{\mbox{\scriptsize rs}}$ 
with generic $x$, 
the second component of 
$(\tau,x)\ast \dot{w}$ does not depend on the 
choice of $\tau$). 
Hence we obtain a birational mapping 
$\ast\dot{w}: \widehat{\frb}\cdots\!\to\widehat{\frb}$,
which gives a commutative diagram
\begin{equation}
\begin{array}{ccccc}
 &\mbox{\small $\mu$} &&\mbox{\small pr}&\\[-2pt]
B\times \frh_{\mbox{\scriptsize reg}} & \iso &
T\times \widehat{\frb}_{\mbox{\scriptsize rs}} & \to&
\widehat{\frb}_{\mbox{\scriptsize rs}} \cr
\mbox{\small$\ast\dot{w}$}\ \vdots\quad & &
\mbox{\small$\ast\dot{w}$}\ \vdots\quad & &
\ \,\vdots \ \mbox{\small$\ast\dot{w}$}\\[-2pt]
\phantom{\mbox{\small$\ast\dot{w}$}}\downarrow\quad & &
\phantom{\mbox{\small$\ast\dot{w}$}}\downarrow\quad & &
\ \downarrow\phantom{\mbox{\small$\ast\dot{w}$}} \\[4pt]
B\times \frh_{\mbox{\scriptsize reg}} & \iso&
T\times \widehat{\frb}_{\mbox{\scriptsize rs}} & \to&
\widehat{\frb}_{\mbox{\scriptsize rs}} 
\end{array}
\end{equation}
of birational mappings. 

For any $b=\tau\exp(\xi)\in B$ with 
generic $\xi\in \widehat{\frn}$, 
let us consider the Gauss decomposition
\begin{equation}
b\,\dot{w}=g(\dot{w},b)^{-1}\,(b\ast\dot{w}),\quad
b\ast\dot{w}=\widetilde{\tau}\exp(\widetilde{\xi})
\quad(\widetilde{\tau}\in T,\ \widetilde{\xi}\in\widehat{\frn}), 
\end{equation}
of $b\,\dot{w}$. 
If we set $x=\Ad(b) h$ so that $\mu(b,h)=(\tau,x)$, 
then we have
\begin{equation}
(\tau,\,x)\ast \dot{w}=\mu((b,h)\ast\dot{w})
=(\widetilde{\tau},\,\widetilde{x}),
\end{equation}
where
\begin{equation}
\widetilde{x}=\Ad(b\ast \dot{w})\Ad(\dot{w}^{-1})h
=\Ad(g(\dot{w},b))\Ad(b)h=\Ad(g(\dot{w},b))x
\end{equation}
since $b\ast\dot{w}=g(\dot{w},b)\,b\,\dot{w}$. 
(As we remarked above, $\widetilde{x}$ does not 
depend on the $T$-component of $b$.)
When $\dot{w}=\ds_i$, in particular, 
from (\ref{eq-phipsi})
we have
\begin{equation}\label{eq-gi}
g_i=\exp\left(\frac{-1}{\tau^{\alpha_i}\psi_i}f_i\right)
=\exp\left(\frac{\alpha_i}{\varphi_i}f_i\right)
=\exp(t_i(x) f_i),
\end{equation} 
hence
\begin{equation}\label{eq-tildex}
\widetilde{x}=\Ad(g_i)x=\exp(t_i(x)\,\ad(f_i))x,\quad
t_i=\frac{\alpha_i}{\varphi_i}. 
\end{equation}
This means that 
the birational action 
$\ast\ds_i : \widehat{\frb}\cdots\!\to\widehat{\frb}$
coincides with $\sigma_i$ of Proposition \ref{prop-sigma}. 
On the other hand, 
since $\widetilde{\tau}=\tau\,\psi_i^{h_i}$, we also have
\begin{equation}\label{eq-tildetau}
\widetilde{\tau}_i=\psi_i\,\tau_i
 =-\dfrac{\varphi_i}{\alpha_i} \tau_i\tau^{-\alpha_i},\qquad 
\widetilde{\tau}_j=\tau_j\ \ (j\ne i).
\end{equation}

In terms of the coordinate rings, these birational action of 
$\dot{w}\in\dot{W}$ 
on $T\times\widehat{\frb}$ and $\widehat{\frb}$
induces the automorphisms $R_{\dot{w}}$
of $\CK(\widehat{\frb})[\tau^{\pm 1}]$
and $\CK(\widehat{\frb})$, respectively. 
Combining Proposition \ref{prop-sigma} with 
(\ref{eq-tildex}) and (\ref{eq-tildetau}), 
we obtain the following theorem. 

\WSCI{\par\medskip}
\begin{thm}\label{thm-R}
The birational action of $\dot{W}$ on $T\times \widehat{\frb}$ is
determined by 
the algebra automorphisms 
$R_{\dot{s}_i}$ $(i\in I)$ of the ring of Laurent polynomials 
$\CK(\widehat{\frb})[\tau^{\pm1}]$ such that
\begin{equation}
\begin{array}{cl}
\smallskip
(\mbox{\rm i}) & 
R_{\dot{s}_i}(\lambda)=\lambda -\alpha_i\br{h_i,\lambda}\quad
(i\in I,\ \lambda \in \frh^\ast)\cr
(\mbox{\rm ii}) & 
R_{\dot{s}_i}(\tau_j)=\tau_j\ \ (j\ne i),\ \ 
R_{\dot{s}_i}(\tau_i)
 =-\dfrac{\varphi_i}{\alpha_i}\,
\tau_i\,{\displaystyle \prod_{k\in I}\tau_k^{-a_{ki}}}
\cr
(\mbox{\rm iii}) & 
R_{s_i}(\psi)=\exp(\dfrac{\lambda_i}{\varphi_i}
\pad(\varphi_i))(\psi)\quad(i\in I,\ \psi\in \CO(\widehat{\frn})),
\end{array}
\end{equation}
where $\lambda_i=\epsilon_i\alpha_i$. 
\end{thm}

\par\medskip
We remark that the birational action of $\dot{W}$ on 
$\widehat{\frb}$ reduces in fact to 
that of the Weyl group $W=W(A)$,  
since the square of $\ast\ds_i=\sigma_i$ gives 
the identity mapping of $\widehat{\frb}$ for each $i\in I$. 
Note, however, that $(\ast\ds_i)^2$ is {\em not} 
the identity as a birational automorphism of 
$T\times \widehat{\frb}$.  
In fact, the action of $R_{\ds_i}$ on the $\tau$-functions  
differ by the factor $-\alpha_i$, 
from that of $s_i$ given in Theorem \ref{thm-B}.

\rem{
In this section, 
we used for the sake of simplicity 
the dual lattice 
$L\subset\Hom_{\BZ}(Q^\vee,\BZ)$ 
of the coroot lattice $Q^\vee$, 
to specify the Cartan subgroup $T$ of $G$. 
A more intrinsic way would be to take an arbitrary 
$\BZ$-lattice of $\frh^{\ast}$ 
such that $Q\subset L$ and $\br{h_i,L}\in\BZ$ for 
all $i\in I$, so that $T=\Spec(\BC[L])$ and 
$\CO(T)=\bigoplus_{\Lambda\in L}\,\BC\, \tau^\Lambda$,
where $\tau^\Lambda$ ($\Lambda\in L$) 
are the formal exponentials. 
Theorems \ref{thm-B} and \ref{thm-R} extends 
naturally to this setting; 
one has only to replace the formulas for $s_i(\tau_j)$ 
and $R_{s_i}(\tau_j)$ by 
\begin{equation}
s_i(\tau^\Lambda)=\varphi_i^{\br{h_i,\Lambda}}\,\tau^{s_i.\Lambda},
\quad
R_{s_i}(\tau^\Lambda)=
\left(-\dfrac{\varphi_i}{\alpha_i}\right)^{\br{h_i,\Lambda}} 
\tau^{s_i.\Lambda}, 
\end{equation}
to obtain a birational realization of $W$ and $\dot{W}$, 
respectively.  
}

\section{Regularity of the $\tau$-cocycle}

In this section, we will prove that the value
$\phi_{w}(\Lambda_j)$ of the $\tau$-cocycle is a polynomial 
in $\CA=\CA_0[\lambda]$ for any $w\in W$ and $j\in I$,
assuming that the generalized Cartan matrix $A$ is 
symmetrizable. 
For this purpose, we may assume that the indexing set $I$ 
is finite, and $\CA_0=\CO(\widehat{\frn})$ as in the previous 
section.   
We will prove the polynomiality of $\phi_{w}(\Lambda_j)$ 
by using the action of the Kac-Moody group $G$ on the 
integrable highest weight $\frg$-modules $L(\Lambda)$
with a dominant integral weight $\Lambda\in P_{+}$. 
Before the proof of the regularity of the $\tau$-cocycle, 
we investigate the difference between the 
birational realization of $\dot{W}$ of Theorem \ref{thm-R}
on $\CK(\widehat{\frb})[\tau^{\pm1}]$ and 
that of $W$ in Theorem \ref{thm-B}.

\subsection{Comparison of the two birational realizations}
In Theorem \ref{thm-R} we gave a realization 
of the group $\dot{W}$ as a group of automorphisms 
of $\CK(\widehat{\frb})[\tau^{\pm1}]$.  
It is also clear that each $R_{\dot{w}}$ preserves the 
subring 
$\CK[\tau^{\pm1}]$ defined by the 
field of fractions $\CK=Q(\CA)$ of $\CA=\BC[\lambda]\otimes\CA_0.$
On this field $\CK$, the automorphisms 
$R_{\dot{s}_i}$ coincide with the automorphism of $s_i : \CK\to\CK$ 
defined by Theorem \ref{thm-A},
while these two actions are different on the $\tau$-functions. 
We describe the difference between the two actions
on the $\tau$-functions.

Note first that 
\begin{equation}
s_i(\tau^\Lambda)=\varphi_i^{\br{h_i,\Lambda}} \tau^{s_i.\Lambda},
\quad
R_{\dot{s}_i^{\pm1}}(\tau^\Lambda)=
\left(\mp \frac{\varphi_i}{\alpha_i}\right)^{\br{h_i,\Lambda}}
\tau^{s_i.\Lambda}
\end{equation}
for each $i\in I$. 
This implies that $w(\tau^\Lambda)$ and $R_{\dot{w}}(\tau^\Lambda)$ 
differ only by a factor which is a rational function of $\alpha_i$ $(i\in I)$,
when $\dot{w}\in\dot{W}$ is a lift of $w\in W$.
Hence, the ratio
$N_{\dot{w}}(\Lambda)={w(\tau^\Lambda)}/{R_{\dot{w}}(\tau^\Lambda)}$
defines a 1-cocycle of $\dot{W}$ with coefficients in 
$\Hom(L,\BQ(\alpha)^\times)$: 
\begin{equation}\begin{array}{l}
\smallskip
N_{\dot{w}}(\Lambda+\Lambda')=N_{\dot{w}}(\Lambda)\,N_{\dot{w}}(\Lambda')
\quad(\Lambda,\Lambda'\in L),\cr
\smallskip
N_1(\Lambda)=1,\quad 
N_{\dot{s}_i}(\Lambda)=(-\alpha_i)^{\br{h_i,\Lambda}}\ \ (i\in I),\cr
N_{\dot{w}\dot{w}'}(\Lambda)= w(N_{\dot{w}'}(\Lambda))\,N_{\dot{w}}(w'.\Lambda)
\quad(\dot{w},\dot{w}'\in \dot{W}). 
\end{array}
\end{equation}
When $\dot{w}=\dot{s}_{j_1}\cdots\dot{s}_{j_p}$, the factors 
$N_{\dot{w}}(\Lambda)$ ($\Lambda\in L$) are determined 
explicitly as 
\begin{equation}
N_{\dot{w}}(\Lambda)
=\prod_{k=1}^p \, (-s_{j_1}\cdots s_{j_{k-1}}(\alpha_{j_k}))
^{\br{h_{j_k},s_{j_{k+1}}\cdots s_{j_{p}}\Lambda}}. 
\end{equation}
Note that, if $w=s_{j_1}\cdots s_{j_p}$ is a reduced decomposition 
of $w\in W$, one has 
\begin{equation}
\Delta_+\cap w(\Delta_-)
=\pbr{\,s_{j_1}\cdots s_{j_{k-1}}(\alpha_{j_k})\ ; k=1,\ldots,p}. 
\end{equation}
Assume furthermore that $\Lambda\in L$ is dominant, i.e., 
$\br{h_i,\Lambda}\ge0$ for each $i\in I$. 
Then one has
\begin{equation}
\br{h_{j_k},s_{j_{k+1}}\cdots s_{j_{p}}\Lambda}\ge 0
\quad(k=1,\ldots,p).
\end{equation}
Hence, $N_{\dot{w}}(\Lambda)$ is a polynomial in $\alpha_{i}$ ($i\in I$).

\WSCI{\par\medskip}
\begin{prop} 
Let $w=s_{j_1}\cdots s_{j_p}$ be a reduced decomposition of 
an element of $w\in W$ and take the corresponding element 
$\dot{w}=\dot{s}_{j_1}\cdots \dot{s}_{j_p}$ in $\dot{W}$. 
Then, for any $\Lambda\in L$, one has
\begin{equation}
R_{\dot{w}}(\tau^\Lambda)=\frac{1}{N_{\dot{w}}(\Lambda)}
w(\tau^\Lambda)=
\frac{\phi_w(\Lambda)}{N_{\dot{w}}(\Lambda)}\,\tau^{w.\Lambda},
\end{equation}
where $N_{\dot{w}}(\Lambda)$ is a rational function in the variables
$\alpha_i$ $(i\in I)$ with coefficients in $\BQ$.  
If $\Lambda\in L$ is dominant, $N_{\dot{w}}(\Lambda)$ is 
a polynomial in $\alpha_i$ $(i\in I)$ with integer coefficients.
\end{prop}

\subsection{Proof of Theorem \ref{thm-C}}

\par\medskip
Recall that the $\BZ$-module of integral weights $P$, and the cone of 
dominant integral weights $P_+$ are defined by
\begin{equation}
\begin{array}{l}
\smallskip
\,P\ =\pbr{\Lambda\in \frh^\ast; \br{h_i,\Lambda}\in \BZ \ \ (i\in I)},
\cr
P_+=\pbr{\Lambda\in \frh^\ast; 
\br{h_i,\Lambda}\in \BZ_{\ge 0} \ (i\in I)},
\end{array}
\end{equation}  
respectively. 
Note that there is a natural surjective  
homomorphism $P\to L$, where $L=\Hom_{\BZ}(Q^\vee,\BZ)$, 
and that, for each $\Lambda\in P$, the values $\tau^\Lambda$ 
depends only on its image in $L$.  
We fix a dominant integral weight $\Lambda\in P_{+}$ and 
consider the expectation value 
\begin{equation}
\br{\Lambda|\cdot|\Lambda }: L(\Lambda)^\circ\times L(\Lambda) \to \BC
\end{equation}
at $\Lambda$.
Here $L(\Lambda)$ is the integrable left $\frg$-module with 
highest weight $\Lambda$, and 
$L(\Lambda)^{\circ}$ is the right $\frg$-module with 
highest weight $\Lambda$,
obtained from $L(\Lambda)$ by 
the anti-involution of $\frg$ such that 
$e_j^\circ=f_j$, $f_j^\circ =e_j$ ($j\in I$) and 
${h}^\circ=h$ for $h\in \frh$. 
Note that 
\begin{equation}
\begin{array}{lll}
\smallskip
\frn|\Lambda\rangle=0,\quad&
h|\Lambda\rangle=|\Lambda\rangle\br{h,\Lambda}\quad&(h\in \frh)\cr
\langle\Lambda|\frn_{-}=0,\quad&
\langle\Lambda|h=\br{h,\Lambda}\langle\Lambda|\quad&(h\in \frh).
\end{array}
\end{equation}
Let $b=\tau\exp(\xi)$ be an element of $B$ with 
$\tau\in T$, $\xi\in \widehat\frn$. 
Then one has 
\begin{equation}
\br{\Lambda|\,b\,|\Lambda}=\br{\Lambda|\tau \exp(\xi)|\Lambda}
=\br{\Lambda|\tau|\Lambda}=\tau^{\Lambda}. 
\end{equation}
Namely, $\tau^\Lambda$ is the expectation value of the Cartan component 
of $b\in B$ at $\Lambda$. 
Let $\dot{w}\in \dot{W}$ and take a generic $b\in B$ such that 
$b\,\dot{w} \in U_{-}\,B$.  
Setting $b\ast \do{w}=\widetilde{\tau}\exp(\widetilde{\xi})$
($\widetilde{\tau}\in T$, $\widetilde{\xi} \in \widehat{\frn}$), 
we have
\begin{equation}
\br{\Lambda|b\ast\dot{w}|\Lambda}=
\br{\Lambda|\,\widetilde{\tau}\exp(\widetilde{\xi})\,|\Lambda}=
\br{\Lambda|\,\widetilde{\tau}\,|\Lambda}
=\widetilde{\tau}^\Lambda=R_{\dot{w}}(\tau^\Lambda),
\end{equation}
where $R_{\dot{w}}(\tau^\Lambda)$ is an 
abbreviation for 
$\prod_{i\in I}(R_{\dot{w}}\tau_i)^{\br{h_i,\Lambda}}$. 
On the other hand, we compute
\begin{equation}
\br{\Lambda|b\ast\dot{w}|\Lambda}
=\br{\Lambda|g(\dot{w},b)\,b\,\dot{w}|\Lambda}
=\br{\Lambda|b\,\dot{w}|\Lambda}
\end{equation}
since $g(\dot{w},b)\in U_{-}$. 
Hence we have
\begin{equation}
R_{\dot{w}}(\tau^\Lambda)=\br{\Lambda|\,b\,\dot{w}|\Lambda}. 
\end{equation}
Note that $w.\Lambda=\Lambda-\beta$, for some $\beta\in Q_{+}$ 
and that there exists an element $F_{\beta}\in U(\frn_-)$ 
of weight $-\beta$ such that $\dot{w}|\Lambda\rangle=F_{\beta} |\Lambda\rangle$.
Hence we have
\begin{equation}
R_{\dot{w}}(\tau^\Lambda)
=\br{\Lambda|\,b\,F_{\beta}|\Lambda}
=\br{\Lambda|\,\tau\exp(\xi)\,F_{\beta}|\Lambda}
=\tau^\Lambda\br{\Lambda|\exp(\xi)_{\beta}\,F_{\beta}|\Lambda},
\end{equation}
where $\exp(\xi)_{\beta}$ stands for the component of weight $\beta$
of $\exp(\xi)$:
\begin{equation}
\exp(\xi)_{\beta}=\sum_{k=0}^\infty \frac{1}{k!} 
\sum_{\beta_1+\cdots+\beta_k=\beta} \xi_{\beta_1}\cdots\xi_{\beta_k}
\qquad(\mbox{finite sum}).
\end{equation}
This implies that $R_{\dot{w}}(\tau^\Lambda)$ is a regular function 
on $B$: 
\begin{equation}
R_{\dot{w}}(\tau^\Lambda) \in \CO(U)\,\tau^\Lambda\subset\CO(B). 
\end{equation}
Hence, under the identification of Proposition \ref{prop-Ad}, we have
\begin{equation}
R_{\dot{w}}(\tau^\Lambda)\in 
\CO(\widehat{\frb}_{\mbox{\scriptsize rs}})[\tau^{\pm1}]
\end{equation}
Since we already know that the left hand side belongs 
to $\CK(\widehat{\frb})\,\tau^{w.\Lambda}$, we obtain
\begin{equation}
R_{\dot{w}}(\tau^\Lambda)\in 
\CO(\widehat{\frb}_{\mbox{\scriptsize rs}})\tau^{w.\Lambda}
\end{equation}
namely, 
\begin{equation}
\frac{\phi_{w}(\Lambda)}{N_{\dot{w}}(\Lambda)}\in 
\CO(\widehat{\frb})[\,\alpha^{-1}\ (\alpha\in \Delta_+)]. 
\end{equation}
Let $w=s_{j_1}\cdots s_{j_p}$ be a reduced decomposition 
of an element $w\in W$ and set $\dot{w}=\dot{s}_{j_1}\cdots \dot{s}_{j_p}$.
In this setting, the normalization factor $N_{\dot{w}}(\Lambda)$ 
is a polynomial in $\alpha_i$ ($i \in I$).  
Hence, we have 
\begin{equation}
\phi_{w}(\Lambda)\in \CO(\widehat{\frb})[\alpha^{-1}\ (\alpha\in \Delta_+)]
\subset \BC(\alpha)\otimes\CO(\widehat{\frn})
=\BC(\alpha)\otimes \CA_0.
\end{equation}
where $\CA_0=\CO(\widehat{\frn})=\Sym(\frn^\ast)$. 

We now prove that 
\begin{equation}
\phi_{w}(\Lambda)\in \CA=
\BC[\alpha]\otimes\CA_0
\qquad(w\in W),
\end{equation}
by the induction on the length $\ell(w)=p$.  
Set $w'=s_{j_2}\cdots s_{j_p}$ so that $w=s_{j_1}w'$.  
Then, by the cocycle property of $\phi$, we have 
\begin{equation}
\phi_{w}(\Lambda)
=s_{j_1}(\phi_{w'}(\Lambda))\,\phi_{s_{j_1}}(w'\Lambda)
=s_{j_1}(\phi_{w'}(\Lambda))\,\varphi_{j_1}^{\br{h_{j_1}, w'\Lambda}}.
\end{equation}
By the induction hypothesis, we have 
$\phi_{w'}(\Lambda)\in \CA$.  
Since $s_{j_1}(\CA)\subset\CA[\varphi_{j_1}^{-1}]$, we see 
$\phi_{w}(\Lambda)\in \CA[\varphi_{j_1}^{-1}]$. 
On the other hand, we already know that 
$\phi_{w}(\Lambda)\in \BC(\alpha)\otimes\CA_0$. 
Since 
$\CA_0=\Sym(\frn^{\ast})$ is a polynomial ring including 
$\varphi_{j_1}$ as an indeterminate, we have
\begin{equation}
\BC[\alpha]\otimes\CA_0[\varphi_{j_1}^{-1}]\ \cap\ 
\BC(\alpha)\otimes \CA_0 =\BC[\alpha]\otimes\CA_0=\CA. 
\end{equation}
Hence $\phi_{w}(\Lambda)\in \CA$, as desired. 
This completes the proof of Theorem \ref{thm-C}. 

\rem{
When the Poisson algebra $\CA_0=\Sym(\frn^{\ast})$ has a $\BZ$-form, 
one can show furthermore that each 
$\phi_w(\Lambda_j)$ is defined over $\BZ$. 
To be more precise, suppose that
there exists a $\BZ$-submodule $\CS\subset\frn^{\ast}$ 
such that \par\smallskip
\begin{tabular}{cl}
(i)&  $\pbr{\CS,\CS}\subset\CS$, and 
$\CS\otimes_{\BZ}\BC=\frn^{\ast}$.
\cr (ii)& $\varphi_j\in \CS$ for any $j\in I$.\cr
(iii)& $\frac{1}{k!}\pad(\varphi_j)^k(\CS)\subset \CS$ 
for any $j\in I$ and $k=0,1,2,\ldots$. 
\end{tabular}
\par\smallskip\noindent 
Then, for any dominant integral weight $\Lambda\in P_{+}$, 
one has  $\phi_{w}(\Lambda)\in\BZ[\lambda][\CS]$ 
for all $w\in W$. 
Namely, each $\phi_{w}(\Lambda)$ is expressed as 
a linear combination 
of products of $\lambda_i=\epsilon_i\alpha_i$ 
($i\in I$) and 
elements of $\CS$, with integer coefficients. 
}

\section{Concluding remarks}

In this paper, we proposed a general method to realize 
an arbitrary Weyl group $W=W(A)$ 
as a group of automorphisms of the 
field of rational functions
$\CK=\BC(\lambda_i, \varphi_i, \pbr{\varphi_i,\varphi_j},\cdots)$, 
starting from a nilpotent Poisson algebra 
generated by $\varphi_i$ ($i\in I$) with the Serre relations. 
Also we discussed how this realization arise 
from the Weyl group of the Kac-Moody group
through the Gauss decomposition. 
\par\medskip
Here we give two remarks about the cases when the 
GCM is of affine type. 
\par\smallskip\noindent
(1) \ When the GCM is of affine type, 
our construction already provides a class 
of discrete integrable systems associated 
with the affine root system, as is discussed 
in our previous paper\WSCI{.\cite{NY} }\PREP{ \cite{NY}.}  
They can be regarded as a universal version containing 
{\em higher} $\varphi$-variables (or $f$-variables \cite{NY}).
\par\smallskip\noindent
(2) \ Birational realizations of affine Weyl groups,
in the sense of this paper, 
arise also as groups of B\"acklund transformations for 
nonlinear differential equations of Painlev\'e type,
obtained by certain reduction from the Drinfeld-Sokolov 
hierarchy of modified type\WSCI{.\cite{DS} }\PREP{ \cite{DS}.} 
\par\medskip
Detail of such specific features of the affine cases will be 
discussed in our forthcoming paper. 

\section*{Acknowledgments}
The authors are grateful to 
Professors T.~Tanisaki, 
Y.~Saito and K.~Iohara for valuable discussions.  
A part of this paper was prepared during the stay 
of one of the authors (M.N.) 
at the Erwin Schr\"odinger International Institute
for Mathematical Physics, 
Vienna, for the research program on Representation 
Theory.  
He would like to express his thanks to Professors V.G.~Kac and A.A.~Kirillov, 
for their hospitality during his stay at ESI. 


\end{document}